\documentclass[11pt]{article}
\usepackage{amsthm}
\usepackage{amsmath}
\usepackage{amssymb}
\usepackage[toc,page]{appendix}

\allowdisplaybreaks

\makeatletter
\@addtoreset{equation}{section}

\makeatother

\makeatletter

\newdimen\mainfontsize \mainfontsize=1\@ptsize pt
\topskip=\mainfontsize
   \@maxdepth=\maxdepth

\makeatother

\theoremstyle{plain}
\newtheorem{thm}{Theorem}[section]
\newtheorem{lem}[thm]{Lemma}

\theoremstyle{definition}
\newtheorem{defn}[thm]{Definition}

\theoremstyle{remark}
\newtheorem{rem}[thm]{Remark}

\allowdisplaybreaks
\numberwithin{thm}{section}
\numberwithin{equation}{section}

\title{ERRATUM: Stochastic evolution equations for large portfolios of stochastic volatility models}
\author{Ben Hambly\footnote{hambly@maths.ox.ac.uk} $\,$ and Nikolaos Kolliopoulos\footnote{kolliopoulos@maths.ox.ac.uk (corresponding author)} \\
Mathematical Institute, University of Oxford}
\date{\today} 

\begin{document}
\maketitle

%
%
%
%
%
%

\begin{abstract}
In the article "Stochastic evolution equations for large portfolios of Stochastic Volatility models" (\cite{HK17}, ArXiv ID: 1701.05640) there is a mistake in the proof of Theorem 3.1. In this erratum we establish a weaker version of this Theorem and then we redevelop the regularity theory for our problem accordingly. This means that most of our regularity results are replaced by slightly weaker ones. We also clarify a point in the proof of a correct result.
\end{abstract}

We will first present the correct results replacing the incorrect ones in a structured way, and then give the proofs. To do this, we require a stronger assumption on the parameters of the CIR volatility process given by $(3.1)$. In the previous version we assumed $\frac{k\theta}{\xi^2}>\frac{3}{4}$. However we now need to impose the stronger condition that $\frac{k\theta}{\xi^{2}} > x^{*}\approx 3.315$, where $x^{*}$ is the largest root of the equation
\begin{equation}
16x^3 - 60x^2 + 24x - 3 = 0,   
\end{equation}
for our results to hold. 
We will also clarify an argument in the (correct) proof of Theorem 4.1 from the original article in the appendix. 

Sections and new results/equations in this erratum will be indexed by numbers preceded by the letter ``E". On the other hand, we will refer to everything else as if we were in the original article.

\section{The corrected main results}
The proof of Theorem~3.1 contains a fatal mistake. We replace the incorrect Theorem~3.1 by the the following:

\begin{thm}\label{thm:3.2}
Suppose that $\sigma_{0}$ is a positive random variable which is bounded away from zero and infinity.
Then $\mathbb{P}$ - almost surely the conditional probability measure $\mathbb{P}(\sigma_{t}\in \cdot\,|\,B_{\cdot}^{0},\,\mathcal{G})$ possesses a continuous density $p_{t}(\cdot\,|\,B_{\cdot}^{0},\,\mathcal{G})$ which is supported in $[0,\infty)$, for all $t > 0$. Moreover, for any $T>0$, any $\alpha\geq0$ and all sufficiently small $q > 1$, we have the following integrability condition
\[
M_{B_{\cdot}^{0},\, \mathcal{G}}^{\alpha}(\cdot) := \sup_{y\geq0}\left(y^{\alpha}p_{\cdot}(y\,|\,B_{\cdot}^{0},\,\mathcal{G})\right) \in L^{q}\left(\Omega \times \left[0, \, T\right]\right)
\]
\end{thm}

The above Theorem has stronger assumptions and gives a weaker result for the volatility density than Theorem~3.1. Thus, the whole regularity theory for $v_{t, C_1}$ needs to be reestablished. We state the results in a slightly different way. The two-dimensional density of the above measure-valued process will belong to the following spaces 
\[
L_{\alpha}=L^{2}\left(\left(\Omega, \, \mathcal{F}, \, \mathbb{P}\right)\times\left[0,\,T\right];\,L_{|y|^{\alpha}}^{2}\left(\mathbb{R}^{+}\times\mathbb{R}\right)\right)
\]
and
\[
H_{\alpha}=L^{2}\left(\left(\Omega, \, \mathcal{F}, \, \mathbb{P}\right)\times\left[0,\,T\right];\,H_{0, w^2(x)}^{1}\left(\mathbb{R}^{+}\right)\times L_{|y|^{\alpha}}^{2}\left(\mathbb{R}\right)\right)
\]
for $\alpha\geq0$ and $w(x) = \min\{1, \, \sqrt{x}\}$ for $x \geq 0$, where we write $L_{g(y)}^{2}$ for the weighted
$L^{2}$ space with weight function $\left\{ g(y):y \in \mathbb{R}\right\}$, and $H_{0, g(x)}^{1}\left(\mathbb{R}^{+}\right)$ for the weighted
$H_{0}^{1}\left(\mathbb{R}^{+}\right)$ space with weight function $\left\{ g(x):x\geq0\right\}$ in the $L^2$ norm of the derivative. Apart from the integrability conditions, a function $u'$ belonging to the second space has to satisfy the boundary condition $\displaystyle{\lim_{x \rightarrow 0^{+}}\left\Vert u'(\cdot, \, x, \, \cdot) \right\Vert_{L_{|y|^{\alpha}}^2\left(\Omega \times \left[0, \, T\right] \times \mathbb{R}\right)} = 0}$. Observe that this definition is not problematic, since $\left\Vert u'(\cdot, \, x, \, \cdot) \right\Vert_{L_{|y|^{\alpha}}^2\left(\Omega \times \left[0, \, T\right] \times \mathbb{R}\right)}$ has to be continuous in $x$ for $x > 0$ (this follows by applying Morrey's inequality \cite{EVAN} away from $x=0$), so changing the value of the above limit gives a different function in an $L^2\left(\Omega \times \mathbb{R}^{+} \times \mathbb{R}^{+} \times \mathbb{R}\right)$ sense. The existence of a density for $v_{t, C_1}$ and its regularity are given in the next Theorem, which replaces Theorem~4.3.

\begin{thm}\label{thm:4.3}
Suppose that $h$ is a continuous function taking values in some compact subset of $\mathbb{R}^{+}$. Suppose also
that given $\mathcal{G}$, $X^{1}_{0}$ has an $L^{2}$-integrable density
$u_{0}(\cdot|\mathcal{G})$ in $\mathbb{R}^{+}$ such that $\mathbb{E}\left[\left \Vert u_{0} \right \Vert_{L^{2}\left(\mathbb{R}^{+}\right)}^{2} \, | \, \mathcal{G} \right]
\in L^{q'}\left(\Omega \right)$ for any $q' > 1$. Suppose finally that $\frac{k_1\theta_1}{\xi_1^2} > x^{*}$ and $\rho_{2,1} \in \left(-1, \, 1\right)$ hold for any possible realization of $C_1 = \left(k_{1},\,\theta_{1},\,\xi_{1},\,r_{1},\,\rho_{1,1},\,
\rho_{2,1}\right)$, and that the random variable $\sigma^{1}_{0}$ is positive and bounded away from zero and infinity.
Then, for any possible realization of $C_{1}$, the measure-valued stochastic process $v_{t,C_{1}}$ has a two-dimensional density $u_{C_{1}}(t,\,\cdot,\,W_{\cdot}^{0},\,B_{\cdot}^{0},\,\mathcal{G})$
belonging to the space $L_{\alpha}$ for any $\alpha \geq 0$. Moreover, when $\rho_3 := \int_{0}^{1}dW^0_tdB^0_t = 0$ and $\mathbb{E}\left[\left \Vert u_{0} \right \Vert_{H_{0, w^2(x)}^{1}\left(\mathbb{R}^{+}\right)}^{2}\, | \, \mathcal{G} \right]
\in L^{q'}\left(\Omega \right)$ for any $q' > 1$, the density belongs to $H_{\alpha}$ as well for any $\alpha \geq 0$.
\end{thm}

Next, we obtain our SPDE exactly as in \cite{HK17}, and we adapt the definition of our initial-boundary value problem to the new regularity results given in the above theorem. For this purpose we define the space $\tilde{L}^2_{\alpha, w} := L_{y^{\alpha}w^2(x)}^{2}\left(\mathbb{R}^{+}
\times\mathbb{R}^{+}\right)$ for any $\alpha \geq 0$, and then we modify Definition 5.1 ($\alpha$-solution to our problem) as follows

\begin{defn}\label{problemdef}
For a given real number $\rho$ and a given value of the coefficient vector $C_{1}$, let $h: \mathbb{R}^{+} \longrightarrow \mathbb{R}^{+}$ be a function having polynomial growth, and $U_0$ be a random function which is extended to be zero outside $\mathbb{R}^{+}$ such that $U_{0}\in L^{2}\left(\Omega;\,\tilde{L}_{\alpha}^{2}
\right)$ and $\left(U_{0}\right)_{x}\in L^{2}\left(\Omega;\,\tilde{L}_{\alpha, w}^{2}
\right)$ for some $\alpha > 0$.
Given $C_{1}$, $\rho$, $\alpha$ and the functions $h$ and $U_{0}$, we say that $u$ is an $\alpha$-solution to our problem when the following are satisfied;

\begin{enumerate}
\item   $u$ is adapted to the filtration $\{\sigma\left(\mathcal{G},\,W_{t}^{0},\,B_{t}^{0}\right):\,t\geq0\}$
and belongs to the space $H_{\alpha}\cap L_{\alpha}$.
\item  $u$ is supported in $\mathbb{R}^{+}$ and satisfies the SPDE
\begin{eqnarray}
u(t,x,y)&=& U_{0}(x,\,y)-r_{1}\int_{0}^{t}\left(u(s,x,y)\right)_{x}ds \nonumber \\
&& \qquad +\frac{1}{2}\int_{0}^{t}h^{2}(y)\left(u(s,x,y)\right)_{x}ds-k_{1}\theta_{1}\int_{0}^{t}\left(u(s,x,y)\right)_{y}ds 
\nonumber \\
&& \qquad +k_{1}\int_{0}^{t}\left(yu(s,\,x,\,y)\right)_{y}ds+\frac{1}{2}\int_{0}^{t}h^{2}(y)\left(u(s,x,y)\right)_{xx}ds
\nonumber \\
& & \qquad +\rho\int_{0}^{t}\left(h\left(y\right)\sqrt{y}u(s,\,x,\,y)\right)_{xy}ds \nonumber \\
& & \qquad +\frac{\xi_{1}^{2}}{2}\int_{0}^{t}\left(yu(s,x,y)\right)_{yy}ds-\rho_{1,1}\int_{0}^{t}h(y)
\left(u(s,x,y)\right)_{x}dW_{s}^{0} \nonumber \\
& & \qquad -\xi_{1}\rho_{2,1}\int_{0}^{t}\left(\sqrt{y}u(s,x,y)\right)_{y}dB_{s}^{0}, \label{eq:5.1}
\end{eqnarray}
for all $x, y \in \mathbb{R}^{+}$, where $u_{y}$, $u_{yy}$
and $u_{xx}$ are considered in the distributional sense over the space of test functions
\[
C_{0}^{test}=\{g\in C_{b}^{2}(\mathbb{R}^{+}\times\mathbb{R}):\,g(0,y)=0, \;\forall y\in\mathbb{R}\}.
\]
\end{enumerate}
\end{defn}

The SPDE of the above definition is satisfied by the density $u_{C_{1}}$ for $\rho = \xi_{1}\rho_{3}\rho_{1,1}\rho_{2,1}$, where $\rho_3$ is the correlation between $W^0$ and $B^0$ (i.e $dW^0_t \cdot dB^0_t = \rho_3dt$), while the regularity properties are also satisfied for all $\alpha > 0$ when $\rho_3 = 0$.

Finally we replace Theorem 5.2, which improves the regularity of our two-dimensional density through the initial-boundary value problem, by the following theorem which differs only in the weighted $L^2$ norm used.  

\begin{thm}\label{thm:5.2}
Fix the value of the coefficient vector $C_1$, the function $h$, the real number $\rho$ and the initial data function $U_{0}$. Let $u$ be an $\alpha$-solution to our problem for all $\alpha\geq 0$. Then, the weak derivative $u_{y}$ of $u$ exists and we have
    \[ 
    u_{y}\in L^{2}\left(\left[0,\,T\right]\times\Omega;\,\tilde{L}^2_{\alpha, w}\right)
    \]
for all $\alpha\geq 2$.
\end{thm}

\section{The main lemmas needed}

To prove Theorem~\ref{thm:3.2} which replaces the incorrect Theorem 3.1 from \cite{HK17}, instead of Lemma 3.5 we need the following stronger result, which contains a generalization of Proposition 2.1.1 from page 78 in \cite{Nualart}.

\begin{lem}\label{lem:3.1}
Let $B$ be a Brownian motion defined on $\left[0, T\right] \times \Omega$ for some $T > 0$ and some probability space $\left(\Omega, \mathcal{F}, \mathbb{P} \right)$, and let $F$ be a random variable which is adapted to the Brownian motion $B$. Suppose also that for some $q, \tilde{q}, r, \tilde{r}, \lambda, \tilde{\lambda} > 1$ and $\alpha \geq 0$, with $q \leq 2$ and $\frac{1}{q} + \frac{1}{\tilde{q}} = \frac{1}{r} + \frac{1}{\tilde{r}} = \frac{q\tilde{r}}{\tilde{\lambda}} + \frac{q\tilde{r}}{\lambda} = 1$, we have $F \in L^{\alpha \lambda \vee \tilde{q}}\left(\Omega\right) \cap \mathbb{D}^{1,2\tilde{\lambda} \vee \frac{q\tilde{r}}{q\tilde{r} - 1}}\left(\Omega\right) \cap \mathbb{D}^{2,qr}\left(\Omega\right)$
and also $\frac{|F|^{m}}{\left\Vert D_{\cdot}F \right\Vert_{L^2\left(\left[0, T\right]\right)}^2} \in L^{q\tilde{r}}\left(\Omega\right)$ for $m \in \{0, \, \alpha\}$. Then, the domain of the adjoint of the derivative operator $D: L^{\tilde{q}}\left(\Omega\right) \cap \mathbb{D}^{1, 2 \vee \frac{q\tilde{r}}{q\tilde{r} - 1}}\left(\Omega \right) \longrightarrow L^{2\vee \frac{q\tilde{r}}{q\tilde{r} - 1}}\left(\Omega ; L^2\left(\left[0, \, T \right]\right)\right)$ (which is an extension of the standard Skorokhood integral $\delta$) contains the process $\frac{D_{\cdot}F}{\left\Vert D_{\cdot}F \right\Vert_{L^2\left(\left[0, T\right]\right)}^2}$. Moreover, $F$ possesses a bounded and continuous density $f_F$ for which we have the estimate
\begin{eqnarray}
&& \sup_{x \in \mathbb{R}}|x|^{\alpha}f_F(x) \nonumber \\
&&\qquad \leq \left \Vert |F|^{\alpha}\delta\left(\frac{D_{\cdot}F}{\left\Vert D_{\cdot}F \right\Vert_{L^2\left(\left[0, T\right]\right)}^2}\right) \right \Vert_{L^{q}\left(\Omega\right)} \nonumber \\
&& \qquad \leq \left(C+2\right)\mathbb{E}^{\frac{1}{qr}}\left[\left\Vert D_{\cdot, \cdot}^2F \right\Vert_{L^2\left(\left[0, T\right]^2\right)}^{qr}\right]\times\mathbb{E}^{\frac{1}{q\tilde{r}}}\left[\left|\frac{|F|^{\alpha}}{\left\Vert D_{\cdot}F \right\Vert_{L^2\left(\left[0, T\right]\right)}^2}\right|^{q\tilde{r}}\right] \nonumber \\
&& \qquad \qquad + C\mathbb{E}^{\frac{1}{qr}}\left[\left\Vert D_{\cdot}F \right\Vert_{L^2\left(\left[0, T\right]\right)}^{qr}\right]\times\mathbb{E}^{\frac{1}{q\tilde{r}}}\left[\left|\frac{|F|^{\alpha}}{\left\Vert D_{\cdot}F \right\Vert_{L^2\left(\left[0, T\right]\right)}^2}\right|^{q\tilde{r}}\right],
\end{eqnarray}
for some $C >0$, with the RHS of the above being finite by our assumptions. 
\end{lem}

Next, to prove Theorem~\ref{thm:4.3} which replaces Theorem 4.3 from \cite{HK17}, the estimate given in Theorem 4.1 is not enough. In particular, we need a stronger estimate for the derivative along with a maximum principle. These are given in the following lemma.

\begin{lem}\label{lem:4dest}
Let $u$ be the density obtained in Theorem~4.1. For some $M > 0$ depending only on $r$ and on some compact interval $\mathcal{I} \subset \mathbb{R}^{+}$ containing the minimum and the maximum of $\sigma_{\cdot}$, we have the estimate
\begin{eqnarray}\label{derivest}
\mathbb{E}\left[ \sup_{0 \leq t \leq T}\left \Vert w(\cdot)u_{x}(t,\,\cdot)\right \Vert_{L^{2}(\mathbb{R}^{+})}^{2}\right] &\leq& Me^{MT}\mathbb{E}\left[ \left \Vert w(\cdot)\left(u_0\right)_{x}(\cdot)\right \Vert_{L^{2}(\mathbb{R}^{+})}^{2}\right] \nonumber \\
&& \quad + Me^{MT}\mathbb{E}\left[ \left \Vert u_0(\cdot)\right \Vert_{L^{2}(\mathbb{R}^{+})}^{2}\right]
\end{eqnarray}
where $w(x) = \min\{1, \, \sqrt{x}\}$ for all $x \geq 0$, provided that the RHS is finite. Then, for some $M' > 0$ depending on $M$ and the initial data, we have the maximum principle
\begin{eqnarray}\label{mp}
\mathbb{E}\left[ \sup_{0 \leq t \leq T}\sup_{x \in \mathbb{R}^{+}}u^2(t,\,x) \right] \leq M'
\end{eqnarray}
\end{lem}

Finally, the proof of Theorem~\ref{thm:5.2} is almost identical to the proof of the corresponding Theorem 5.2 in \cite{HK17}. The only difference is that the $\delta$-identity contains two extra non-derivative terms, while the weight $w^2(x)$ is introduced to the norms and inner products involved in all the other terms. This is not a problem because the two extra terms are
\begin{eqnarray}
\int_{0}^{t}\left\Vert \mathbb{I}_{\left[0, \, 1\right] \times \mathbb{R}}(\cdot)I_{\epsilon,1}(s,\cdot)\right\Vert _{L^2\left(\Omega^0 ; \, \tilde{L}_{\delta}^2\right)}^2ds
\end{eqnarray}
and
\begin{eqnarray}
\int_{0}^{t}\left\langle \frac{\partial}{\partial x}I_{\epsilon,h^{2}(z)}(s,\,\cdot),\,\mathbb{I}_{\left[0, \, 1\right] \times \mathbb{R}}(\cdot)I_{\epsilon,1}(s,\,\cdot)\right\rangle _{L^{2}\left(\Omega^0;\,\tilde{L}_{\delta}^2\right)}ds
\end{eqnarray}
which do not explode as $\epsilon \rightarrow 0^{+}$ (by Lemma 5.3 and our regularity assumptions), while for all the other terms we use Lemmas 5.3 and 5.4 for a slightly differently weighted measure $\mu$ which gives the weight $w^2(x)$ to the norms and inner products involved. Therefore, we only need to prove the corrected $\delta$-identity which is stated below.

\begin{lem}[\textbf{the $\delta$-identity}]

The following identity holds for any $\delta > 1$
\begin{eqnarray}
\left\Vert I_{\epsilon,1}(t,\,\cdot)\right\Vert _{L^{2}\left(\Omega^0;\,\tilde{L}_{\delta, w}^2\right)}^{2} &=& \left\Vert \int_{\mathcal{D}}U_{0}(\cdot,\,z)\phi_{\epsilon}(z,\,\cdot)dz\right\Vert _{L^{2}\left(\Omega^0;\,\tilde{L}_{\delta, w}^2\right)}^{2} \nonumber \\
& & \qquad +r_{1}\int_{0}^{t}\left\Vert \mathbb{I}_{\left[0, \, 1\right] \times \mathbb{R}}(\cdot)I_{\epsilon,1}(s,\cdot)\right\Vert _{L^2\left(\Omega^0 ; \, \tilde{L}_{\delta}^2\right)}^2ds \nonumber \\
& & \qquad +\int_{0}^{t}\left\langle \frac{\partial}{\partial x}I_{\epsilon,h^{2}(z)}(s,\,\cdot),\,I_{\epsilon,1}(s,\,\cdot)\right\rangle _{L^{2}\left(\Omega^0;\,\tilde{L}_{\delta, w}^2\right)}ds \nonumber \\
& & \qquad + \delta \left( k_{1}\theta_{1} - \frac{\xi_1^2}{4}\right)\int_{0}^{t}\left\langle I_{\epsilon, z^{-\frac{1}{2}}}(s,\,\cdot),\,I_{\epsilon,1}(s,\,\cdot)\right\rangle _{L^{2}\left(\Omega^0;\,\tilde{L}^2_{\delta-1, w}\right)}ds \nonumber \\
& & \qquad +\left(k_{1}\theta_{1} - \frac{\xi_1^2}{4}\right)\int_{0}^{t}\left\langle I_{\epsilon, z^{-\frac{1}{2}}}(s,\,\cdot),\,\frac{\partial}{\partial y}I_{\epsilon,1}(s,\,\cdot)\right\rangle _{L^{2}\left(\Omega^0;\,\tilde{L}_{\delta, w}^2\right)}ds \nonumber \\
& & \qquad -\delta k_{1}\int_{0}^{t}\left\langle I_{\epsilon, z^{\frac{1}{2}}}(s,\,\cdot),\,I_{\epsilon,1}(s,\,\cdot)\right\rangle _{L^{2}\left(\Omega^0;\,\tilde{L}^2_{\delta-1, w}\right)}ds \nonumber \\
& & \qquad -k_{1}\int_{0}^{t}\left\langle I_{\epsilon,z^{\frac{1}{2}}}(s,\,\cdot),\,\frac{\partial}{\partial y}I_{\epsilon,1}(s,\,\cdot)\right\rangle _{L^{2}\left(\Omega^0;\,\tilde{L}_{\delta, w}^2\right)}ds \nonumber \\
& & \qquad -\int_{0}^{t}\left\langle \frac{\partial}{\partial x}I_{\epsilon,h^{2}(z)}(s,\,\cdot),\,\frac{\partial}{\partial x}I_{\epsilon,1}(s,\,\cdot)\right\rangle _{L^{2}\left(\Omega^0;\,\tilde{L}_{\delta, w}^2\right)}ds \nonumber \\
& & \qquad -\int_{0}^{t}\left\langle \frac{\partial}{\partial x}I_{\epsilon,h^{2}(z)}(s,\,\cdot),\,\mathbb{I}_{\left[0, \, 1\right] \times \mathbb{R}}(\cdot)I_{\epsilon,1}(s,\,\cdot)\right\rangle _{L^{2}\left(\Omega^0;\,\tilde{L}_{\delta}^2\right)}ds \nonumber \\
& & \qquad -\delta\rho\int_{0}^{t}\left\langle \frac{\partial}{\partial x}I_{\epsilon,h\left(z\right)}(s,\cdot),\,I_{\epsilon,1}(s,\cdot)\right\rangle _{L^2\left(\Omega^0;\, \tilde{L}_{\delta-1, w}^2 \right)}ds \nonumber \\
& & \qquad +\rho_{1,1}^{2}\int_{0}^{t}\left\Vert \frac{\partial}{\partial x}I_{\epsilon,h(z)}(s,\,\cdot)\right\Vert _{L^{2}\left(\Omega^0;\,\tilde{L}_{\delta, w}^2\right)}^{2}ds
\nonumber \\
& & \qquad +\delta(\delta-1)\frac{\xi_{1}^{2}}{8}\int_{0}^{t}\left\Vert I_{\epsilon,1}(s,\,\cdot)\right\Vert_{L^2\left(\Omega^0;\,\tilde{L}^2_{\delta-2, w}\right)}^{2}ds
\nonumber \\
& & \qquad -\frac{\xi_{1}^{2}}{4}\left(1-\rho_{2,1}^{2}\right)\int_{0}^{t}\left\Vert \frac{\partial}{\partial y}I_{\epsilon,1}(s,\,\cdot)\right\Vert _{L^{2}\left(\Omega^0;\,\tilde{L}_{\delta, w}^2\right)}^{2}ds. \nonumber \\
& & \qquad - \left(\rho - \xi_{1}\rho_{3}\rho_{1,1}\rho_{2,1}\right) \nonumber \\
& & \qquad \qquad \times\int_{0}^{t}\left\langle \frac{\partial}{\partial x}I_{\epsilon,h\left(z\right)}(s,\cdot),\,\frac{\partial}{\partial y}I_{\epsilon,1}(s,\cdot)\right\rangle _{L^2\left(\Omega^0; \, \tilde{L}_{\delta, w}^2\right)}ds. \label{eq:5.17}
\end{eqnarray} 
All the terms in the above identity are finite.
\end{lem}

\section{Proofs}

\begin{proof}[\textbf{Proof of Lemma E2.1}]
Let $\{F_n: n \in \mathbb{N}\}$ be a sequence of regular enough (in terms of Malliavin differentiability) random variables which are adapted to the Brownian motion $B_{\cdot}$ in $\left[0, \, T\right]$, with $F_n \longrightarrow F$ in $L^{\alpha \lambda \vee \tilde{q}}\left(\Omega\right) \cap \mathbb{D}^{1,2\tilde{\lambda} \vee \frac{q\tilde{r}}{q\tilde{r} - 1}}\left(\Omega\right) \cap \mathbb{D}^{2,qr}\left(\Omega\right)$. Then, $\frac{D_{\cdot}F_n}{\left\Vert D_{\cdot}F_n \right\Vert_{L^2\left(\left[0, T\right]\right)}^2 + \epsilon}$ belongs to the domain of the standard Skorokhod integral $\delta$ for any $\epsilon > 0$, and for any $m \leq \alpha$, by a well-known property of $\delta$ (see property (4) on page 40 in \cite{Nualart}) we have the following relationship, 
\begin{eqnarray}\label{newlemma1}
&&|F|^m\delta\left(\frac{D_{\cdot}F_n}{\left\Vert D_{\cdot}F_n \right\Vert_{L^2\left(\left[0, T\right]\right)}^2 + \epsilon}\right) \nonumber \\
&& \qquad = \frac{|F|^m\delta\left(D_{\cdot}F_n\right)}{\left\Vert D_{\cdot}F_n \right\Vert_{L^2\left(\left[0, T\right]\right)}^2 + \epsilon} + |F|^m\int_{0}^{T}D_{s}F_nD_{s}\left(\frac{1}{\left\Vert D_{\cdot}F_n \right\Vert_{L^2\left(\left[0, T\right]\right)}^2 + \epsilon}\right)ds \nonumber \\
&& \qquad = \frac{|F|^m\delta\left(D_{\cdot}F_n\right)}{\left\Vert D_{\cdot}F_n \right\Vert_{L^2\left(\left[0, T\right]\right)}^2 + \epsilon} + |F|^m\int_{0}^{T}D_{s}F_n\frac{-\int_{0}^{T}2D_{s'}F_n \cdot D_{s',s}^2F_nds'}{\left(\left\Vert D_{\cdot}F_n \right\Vert_{L^2\left(\left[0, T\right]\right)}^2 + \epsilon\right)^2}ds \nonumber \\
&& \qquad = \frac{|F|^m\delta\left(D_{\cdot}F_n\right)}{\left\Vert D_{\cdot}F_n \right\Vert_{L^2\left(\left[0, T\right]\right)}^2 + \epsilon} - 2|F|^m\frac{\int_{0}^{T}\int_{0}^{T}D_{s}F_n \cdot D_{s'}F_n \cdot D_{s',s}^2F_nds'ds}{\left(\left\Vert D_{\cdot}F_n \right\Vert_{L^2\left(\left[0, T\right]\right)}^2 + \epsilon\right)^2}. \nonumber \\
\end{eqnarray} 
Thus, by the triangle inequality, a boundedness property of the operator $\delta$ (see Proposition 1.5.4 on page 69 in \cite{Nualart}) and H\"older's inequality, we have that
\begin{eqnarray}\label{newlemma2}
&&\left\Vert |F|^m\delta\left(\frac{D_{\cdot}F_n}{\left\Vert D_{\cdot}F_n \right\Vert_{L^2\left(\left[0, T\right]\right)}^2 + \epsilon}\right)\right\Vert_{L^q\left(\Omega\right)} \nonumber \\
&& \qquad \leq \mathbb{E}^{\frac{1}{q}}\left[\left|\frac{|F|^m\delta\left(D_{\cdot}F_n\right)}{\left\Vert D_{\cdot}F_n \right\Vert_{L^2\left(\left[0, T\right]\right)}^2 + \epsilon}\right|^q\right] \nonumber \\
&& \qquad \qquad + 2\mathbb{E}^{\frac{1}{q}}\left[\left||F|^m\frac{\int_{0}^{T}\int_{0}^{T}D_{s}F_n \cdot D_{s'}F_n \cdot D_{s',s}^2F_nds'ds}{\left(\left\Vert D_{\cdot}F_n \right\Vert_{L^2\left(\left[0, T\right]\right)}^2 + \epsilon\right)^2}\right|^q\right] \nonumber \\
&& \qquad \leq \mathbb{E}^{\frac{1}{qr}}\left[\left|\delta\left(D_{\cdot}F_n\right)\right|^{qr}\right]\times\mathbb{E}^{\frac{1}{q\tilde{r}}}\left[\left|\frac{|F|^m}{\left\Vert D_{\cdot}F_n \right\Vert_{L^2\left(\left[0, T\right]\right)}^2 + \epsilon}\right|^{q\tilde{r}}\right] \nonumber \\
&& \qquad \qquad + 2\mathbb{E}^{\frac{1}{q}}\left[\left||F|^m\frac{\int_{0}^{T}\int_{0}^{T}D_{s}F_n \cdot D_{s'}F_n \cdot D_{s',s}^2F_nds'ds}{\left(\left\Vert D_{\cdot}F_n \right\Vert_{L^2\left(\left[0, T\right]\right)}^2 + \epsilon\right)^2}\right|^q\right] \nonumber \\
&& \qquad \leq C_{qr}\mathbb{E}^{\frac{1}{qr}}\left[\left(\int_{0}^{T}\int_{0}^{T}\left|D_{s', s}^2F_n\right|^{2}ds'ds\right)^{\frac{qr}{2}}\right]\times\mathbb{E}^{\frac{1}{q\tilde{r}}}\left[\left|\frac{|F|^m}{\left\Vert D_{\cdot}F_n \right\Vert_{L^2\left(\left[0, T\right]\right)}^2 + \epsilon}\right|^{q\tilde{r}}\right] \nonumber \\
&& \qquad \qquad + C_{qr}\mathbb{E}^{\frac{1}{qr}}\left[\left(\int_{0}^{T}\left|D_{s'}F_n\right|^{2}ds'\right)^{\frac{qr}{2}}\right]\times\mathbb{E}^{\frac{1}{q\tilde{r}}}\left[\left|\frac{|F|^m}{\left\Vert D_{\cdot}F_n \right\Vert_{L^2\left(\left[0, T\right]\right)}^2 + \epsilon}\right|^{q\tilde{r}}\right] \nonumber \\
&& \qquad \qquad+ 2\mathbb{E}^{\frac{1}{q}}\left[|F|^m\left(\frac{\left\Vert D_{\cdot}F_n \right\Vert_{L^2\left(\left[0, T\right]\right)}^2\left(\int_{0}^{T}\int_{0}^{T}\left(D_{s',s}^2F_n\right)^2ds'ds\right)^{\frac{1}{2}}}{\left(\left\Vert D_{\cdot}F_n \right\Vert_{L^2\left(\left[0, T\right]\right)}^2 + \epsilon\right)^2}\right)^q\right] \nonumber \\
&& \qquad \leq C_{qr}\mathbb{E}^{\frac{1}{qr}}\left[\left(\int_{0}^{T}\int_{0}^{T}\left|D_{s', s}^2F_n\right|^{2}ds'ds\right)^{\frac{qr}{2}}\right]\times\mathbb{E}^{\frac{1}{q\tilde{r}}}\left[\left|\frac{|F|^m}{\left\Vert D_{\cdot}F_n \right\Vert_{L^2\left(\left[0, T\right]\right)}^2 + \epsilon}\right|^{q\tilde{r}}\right] \nonumber \\
&& \qquad \qquad + C_{qr}\mathbb{E}^{\frac{1}{qr}}\left[\left(\int_{0}^{T}\left|D_{s'}F_n\right|^{2}ds'\right)^{\frac{qr}{2}}\right]\times\mathbb{E}^{\frac{1}{q\tilde{r}}}\left[\left|\frac{|F|^m}{\left\Vert D_{\cdot}F_n \right\Vert_{L^2\left(\left[0, T\right]\right)}^2 + \epsilon}\right|^{q\tilde{r}}\right] \nonumber \\
&& \qquad \qquad + 2\mathbb{E}^{\frac{1}{q}}\left[|F|^m\left(\frac{\left(\int_{0}^{T}\int_{0}^{T}\left(D_{s',s}^2F_n\right)^2ds'ds\right)^{\frac{1}{2}}}{\left\Vert D_{\cdot}F_n \right\Vert_{L^2\left(\left[0, T\right]\right)}^2 + \epsilon}\right)^q\right] \nonumber \\
&& \qquad \leq \left(C_{qr} + 2\right)\mathbb{E}^{\frac{1}{qr}}\left[\left(\int_{0}^{T}\int_{0}^{T}\left|D_{s', s}^2F_n\right|^{2}ds'ds\right)^{\frac{qr}{2}}\right] \times\mathbb{E}^{\frac{1}{q\tilde{r}}}\left[\left|\frac{|F|^m}{\left\Vert D_{\cdot}F_n \right\Vert_{L^2\left(\left[0, T\right]\right)}^2 + \epsilon}\right|^{q\tilde{r}}\right] \nonumber \\
&& \qquad \qquad + C_{qr}\mathbb{E}^{\frac{1}{qr}}\left[\left(\int_{0}^{T}\left|D_{s'}F_n\right|^{2}ds'\right)^{\frac{qr}{2}}\right]\times\mathbb{E}^{\frac{1}{q\tilde{r}}}\left[\left|\frac{|F|^m}{\left\Vert D_{\cdot}F_n \right\Vert_{L^2\left(\left[0, T\right]\right)}^2 + \epsilon}\right|^{q\tilde{r}}\right] \nonumber \\
&& \qquad = \left(C_{qr}+2\right)\mathbb{E}^{\frac{1}{qr}}\left[\left\Vert D_{\cdot, \cdot}^2F_n \right\Vert_{L^2\left(\left[0, T\right]^2\right)}^{qr}\right]\times\mathbb{E}^{\frac{1}{q\tilde{r}}}\left[\left|\frac{|F|^m}{\left\Vert D_{\cdot}F_n \right\Vert_{L^2\left(\left[0, T\right]\right)}^2 + \epsilon}\right|^{q\tilde{r}}\right] \nonumber \\
&& \qquad \qquad + C_{qr}\mathbb{E}^{\frac{1}{qr}}\left[\left\Vert D_{\cdot}F_n \right\Vert_{L^2\left(\left[0, T\right]\right)}^{qr}\right]\times\mathbb{E}^{\frac{1}{q\tilde{r}}}\left[\left|\frac{|F|^m}{\left\Vert D_{\cdot}F_n \right\Vert_{L^2\left(\left[0, T\right]\right)}^2 + \epsilon}\right|^{q\tilde{r}}\right], \nonumber \\
\end{eqnarray}
for $r,\tilde{r}>1$ such that $\frac{1}{r}+\frac{1}{\tilde{r}}=1$.
Then, for a fixed $\epsilon > 0$, we can use the Lipschitz continuity of $\frac{1}{\epsilon + x^2}$, H\"older's inequality and our assumptions, to show that the last expression converges as $n \longrightarrow +\infty$ to the finite quantity
\begin{eqnarray}
&&\left(C_{qr}+2\right)\mathbb{E}^{\frac{1}{qr}}\left[\left\Vert D_{\cdot, \cdot}^2F \right\Vert_{L^2\left(\left[0, T\right]^2\right)}^{qr}\right]\times\mathbb{E}^{\frac{1}{q\tilde{r}}}\left[\left|\frac{|F|^m}{\left\Vert D_{\cdot}F \right\Vert_{L^2\left(\left[0, T\right]\right)}^2 + \epsilon}\right|^{q\tilde{r}}\right] \nonumber \\
&& \qquad + C_{qr}\mathbb{E}^{\frac{1}{qr}}\left[\left\Vert D_{\cdot}F \right\Vert_{L^2\left(\left[0, T\right]\right)}^{qr}\right]\times\mathbb{E}^{\frac{1}{q\tilde{r}}}\left[\left|\frac{|F|^m}{\left\Vert D_{\cdot}F \right\Vert_{L^2\left(\left[0, T\right]\right)}^2 + \epsilon}\right|^{q\tilde{r}}\right], \nonumber \\
\end{eqnarray}
which implies that for a sequence $\{k_n: n \in \mathbb{N}\} \subset \mathbb{N}$ we have also  
\begin{eqnarray}
|F|^m\delta\left(\frac{D_{\cdot}F_{k_n}}{\left\Vert D_{\cdot}F_{k_n} \right\Vert_{L^2\left(\left[0, T\right]\right)}^2 + \epsilon}\right) \longrightarrow \delta_{F, \epsilon}^{m} 
\end{eqnarray}
weakly in $L^{q}\left(\Omega\right)$ as $n \longrightarrow +\infty$, for some $\delta_{F, \epsilon}^{m}$. Moreover, when $m = 0$, for any $\tilde{F} \in L^{\tilde{q}}\left(\Omega\right) \cap \mathbb{D}^{1, 2 \vee \frac{q\tilde{r}}{q\tilde{r} - 1}}\left(\Omega \right)$ we have
\begin{eqnarray}
\mathbb{E}\left[ \tilde{F} \delta_{F, \epsilon}^{0} \right] &=& \lim_{n \rightarrow + \infty} \mathbb{E}\left[ \tilde{F} \delta\left(\frac{D_{\cdot}F_{k_n}}{\left\Vert D_{\cdot}F_{k_n} \right\Vert_{L^2\left(\left[0, T\right]\right)}^2 + \epsilon}\right) \right] \nonumber \\
&=& \lim_{n \rightarrow + \infty} \mathbb{E}\left[ \int_{0}^{T} D_s\tilde{F} \frac{D_{s}F_{k_n}}{\left\Vert D_{\cdot}F_{k_n} \right\Vert_{L^2\left(\left[0, T\right]\right)}^2 + \epsilon}ds \right] \nonumber \\
&=& \mathbb{E}\left[ \int_{0}^{T} D_s\tilde{F} \frac{D_{s}F}{\left\Vert D_{\cdot}F \right\Vert_{L^2\left(\left[0, T\right]\right)}^2 + \epsilon}ds \right]. \nonumber \\
\end{eqnarray}
To see this we observe that
\begin{eqnarray}
&& \left|\mathbb{E}\left[ \int_{0}^{T} D_s\tilde{F} \left(\frac{D_{s}F_{k_n}}{\left\Vert D_{\cdot}F_{k_n} \right\Vert_{L^2\left(\left[0, T\right]\right)}^2 + \epsilon} - \frac{D_{s}F}{\left\Vert D_{\cdot}F \right\Vert_{L^2\left(\left[0, T\right]\right)}^2 + \epsilon}\right)ds \right]\right| \nonumber \\
&& \qquad \leq \mathbb{E}\left[ \left(\frac{1}{\left\Vert D_{\cdot}F_{k_n} \right\Vert_{L^2\left(\left[0, T\right]\right)}^2 + \epsilon}\right)\int_{0}^{T} \left|D_s\tilde{F}\right| \left| D_sF - D_sF_{k_n} \right|ds \right] \nonumber \\
&& \qquad \qquad + \mathbb{E}\left[ \left|\frac{1}{\left\Vert D_{\cdot}F_{k_n} \right\Vert_{L^2\left(\left[0, T\right]\right)}^2 + \epsilon} - \frac{1}{\left\Vert D_{\cdot}F \right\Vert_{L^2\left(\left[0, T\right]\right)}^2 + \epsilon}\right|\int_{0}^{T} D_s\tilde{F}D_{s}F ds \right] \nonumber \\
&& \qquad \leq \frac{1}{\epsilon}\mathbb{E}^{1 - \frac{1}{q\tilde{r}}}\left[\left\Vert D_{\cdot}\tilde{F} \right\Vert_{L^2\left(\left[0, T\right]\right)}^{\frac{q\tilde{r}}{q\tilde{r} - 1}} \right]\mathbb{E}^{\frac{1}{q\tilde{r}}}\left[ \left\Vert D_sF - D_sF_{k_n} \right\Vert_{L^2\left(\left[0, T\right]\right)}^{q\tilde{r}} \right] \nonumber \\
&& \qquad \qquad + \mathbb{E}\left[ \frac{\left|\left\Vert D_{\cdot}F_{k_n} \right\Vert_{L^2\left(\left[0, T\right]\right)}^2 - \left\Vert D_{\cdot}F \right\Vert_{L^2\left(\left[0, T\right]\right)}^2\right|\left\Vert D_{\cdot}\tilde{F} \right\Vert_{L^2\left(\left[0, T\right]\right)}\left\Vert D_{\cdot}F \right\Vert_{L^2\left(\left[0, T\right]\right)}}{\left(\left\Vert D_{\cdot}F_{k_n} \right\Vert_{L^2\left(\left[0, T\right]\right)}^2 + \epsilon\right)\left(\left\Vert D_{\cdot}F \right\Vert_{L^2\left(\left[0, T\right]\right)}^2 + \epsilon\right)} \right] \nonumber \\
&& \qquad \leq \frac{1}{\epsilon}\mathbb{E}^{1 - \frac{1}{q\tilde{r}}}\left[\left\Vert D_{\cdot}\tilde{F} \right\Vert_{L^2\left(\left[0, T\right]\right)}^{\frac{q\tilde{r}}{q\tilde{r} - 1}} \right]\mathbb{E}^{\frac{1}{q\tilde{r}}}\left[ \left\Vert D_sF - D_sF_{k_n} \right\Vert_{L^2\left(\left[0, T\right]\right)}^{q\tilde{r}} \right] \nonumber \\
&& \qquad \qquad + C_{\epsilon}\mathbb{E}\left[\left|\left\Vert D_{\cdot}F_{k_n} \right\Vert_{L^2\left(\left[0, T\right]\right)} - \left\Vert D_{\cdot}F \right\Vert_{L^2\left(\left[0, T\right]\right)}\right|\left\Vert D_{\cdot}\tilde{F} \right\Vert_{L^2\left(\left[0, T\right]\right)} \right] \nonumber \\
&& \qquad \leq \left(\frac{1}{\epsilon} + C_{\epsilon}\right)\mathbb{E}^{1 - \frac{1}{q\tilde{r}}}\left[\left\Vert D_{\cdot}\tilde{F} \right\Vert_{L^2\left(\left[0, T\right]\right)}^{\frac{q\tilde{r}}{q\tilde{r} - 1}} \right]\mathbb{E}^{\frac{1}{q\tilde{r}}}\left[ \left\Vert D_sF - D_sF_{k_n} \right\Vert_{L^2\left(\left[0, T\right]\right)}^{q\tilde{r}} \right] \nonumber \\
\end{eqnarray}
for some $C_{\epsilon} > 0$, which converges to zero as $n\to\infty$. Furthermore, by a density argument, we can easily show that there is a unique weak limit $\delta_{F, \epsilon}^{0}$ and, since the subsequence $\{k_n: n \in \mathbb{N}\}$ can be taken to be the same for both $m = 0$ and $m = \alpha$, we can also show that $\delta_{F, \epsilon}^{\alpha} = |F|^{\alpha}\delta_{F, \epsilon}^{0}$. Hence, we can define $\delta\left(\frac{D_{\cdot}F}{\left\Vert D_{\cdot}F \right\Vert_{L^2\left(\left[0, T\right]\right)}^2 + \epsilon}\right) := \delta_{F, \epsilon}^{0}$ and then have $\delta_{F, \epsilon}^{\alpha} = |F|^{\alpha}\delta\left(\frac{D_{\cdot}F}{\left\Vert D_{\cdot}F \right\Vert_{L^2\left(\left[0, T\right]\right)}^2 + \epsilon}\right)$ as well. Thus, we can also use Fatou's lemma to estimate
\begin{eqnarray}
&&\left \Vert |F|^m\delta\left(\frac{D_{\cdot}F}{\left\Vert D_{\cdot}F \right\Vert_{L^2\left(\left[0, T\right]\right)}^2 + \epsilon}\right) \right \Vert_{L^q\left(\Omega\right)} \nonumber \\
&& \qquad \leq \liminf_{n \rightarrow +\infty}\left \Vert |F|^m\delta\left(\frac{D_{\cdot}F_{k_n}}{\left\Vert D_{\cdot}F_{k_n} \right\Vert_{L^2\left(\left[0, T\right]\right)}^2 + \epsilon}\right) \right \Vert_{L^q\left(\Omega\right)} \nonumber \\
&& \qquad \leq \left(C_{qr}+2\right)\liminf_{n \rightarrow +\infty}\mathbb{E}^{\frac{1}{qr}}\left[\left\Vert D_{\cdot, \cdot}^2F_{k_n} \right\Vert_{L^2\left(\left[0, T\right]^2\right)}^{qr}\right] \nonumber \\
&& \qquad \qquad \qquad \times\mathbb{E}^{\frac{1}{q\tilde{r}}}\left[\left|\frac{|F|^m}{\left\Vert D_{\cdot}F_{k_n} \right\Vert_{L^2\left(\left[0, T\right]\right)}^2 + \epsilon}\right|^{q\tilde{r}}\right] \nonumber \\
&& \qquad \qquad + C_{qr}\liminf_{n \rightarrow +\infty}\mathbb{E}^{\frac{1}{qr}}\left[\left\Vert D_{\cdot}F_{k_n} \right\Vert_{L^2\left(\left[0, T\right]\right)}^{qr}\right]\times\mathbb{E}^{\frac{1}{q\tilde{r}}}\left[\left|\frac{|F|^m}{\left\Vert D_{\cdot}F_{k_n} \right\Vert_{L^2\left(\left[0, T\right]\right)}^2 + \epsilon}\right|^{q\tilde{r}}\right] \nonumber \\
&& \qquad = \left(C_{qr}+2\right)\mathbb{E}^{\frac{1}{qr}}\left[\left\Vert D_{\cdot, \cdot}^2F \right\Vert_{L^2\left(\left[0, T\right]^2\right)}^{qr}\right]\times\mathbb{E}^{\frac{1}{q\tilde{r}}}\left[\left|\frac{|F|^m}{\left\Vert D_{\cdot}F \right\Vert_{L^2\left(\left[0, T\right]\right)}^2 + \epsilon}\right|^{q\tilde{r}}\right] \nonumber \\
&& \qquad \qquad + C_{qr}\mathbb{E}^{\frac{1}{qr}}\left[\left\Vert D_{\cdot}F \right\Vert_{L^2\left(\left[0, T\right]\right)}^{qr}\right]\times\mathbb{E}^{\frac{1}{q\tilde{r}}}\left[\left|\frac{|F|^m}{\left\Vert D_{\cdot}F \right\Vert_{L^2\left(\left[0, T\right]\right)}^2 + \epsilon}\right|^{q\tilde{r}}\right] \nonumber \\
\end{eqnarray}
for both $m = 0$ and $m = \alpha$. This means that we can take $\epsilon \downarrow 0$ and repeat the previous argument (where this time, we use the Monotone Convergence Theorem to compute the limits) to deduce that $\delta\left(\frac{D_{\cdot}F}{\left\Vert D_{\cdot}F \right\Vert_{L^2\left(\left[0, T\right]\right)}^2}\right)$ can be defined such that  
\begin{eqnarray}\label{property}
\mathbb{E}\left[ \tilde{F} \delta\left(\frac{D_{\cdot}F}{\left\Vert D_{\cdot}F \right\Vert_{L^2\left(\left[0, T\right]\right)}^2}\right) \right] = \mathbb{E}\left[ \int_{0}^{T} D_t\tilde{F} \frac{D_{t}F}{\left\Vert D_{\cdot}F \right\Vert_{L^2\left(\left[0, T\right]\right)}^2}dt \right]
\end{eqnarray}
for any $\tilde{F} \in L^{\tilde{q}}\left(\Omega\right) \cap \mathbb{D}^{1, 2 \vee \frac{q\tilde{r}}{q\tilde{r} - 1}}\left(\Omega \right)$ and
\begin{eqnarray}\label{densityest1}
&&\left \Vert |F|^m\delta\left(\frac{D_{\cdot}F}{\left\Vert D_{\cdot}F \right\Vert_{L^2\left(\left[0, T\right]\right)}^2}\right) \right \Vert_{L^q\left(\Omega\right)} \nonumber \\
&& \qquad \leq \left(C_{qr}+2\right)\mathbb{E}^{\frac{1}{qr}}\left[\left\Vert D_{\cdot, \cdot}^2F \right\Vert_{L^2\left(\left[0, T\right]^2\right)}^{qr}\right]\times\mathbb{E}^{\frac{1}{q\tilde{r}}}\left[\left|\frac{|F|^m}{\left\Vert D_{\cdot}F \right\Vert_{L^2\left(\left[0, T\right]\right)}^2}\right|^{q\tilde{r}}\right] \nonumber \\
&& \qquad \qquad + C_{qr}\mathbb{E}^{\frac{1}{qr}}\left[\left\Vert D_{\cdot}F \right\Vert_{L^2\left(\left[0, T\right]\right)}^{qr}\right]\times\mathbb{E}^{\frac{1}{q\tilde{r}}}\left[\left|\frac{|F|^m}{\left\Vert D_{\cdot}F \right\Vert_{L^2\left(\left[0, T\right]\right)}^2}\right|^{q\tilde{r}}\right]. \nonumber \\
\end{eqnarray}
for both $m = 0$ and $m = \alpha$. The finiteness of the RHS in both \eqref{property} and \eqref{densityest1}, which allows us to obtain these relations by using the Monotone Convergence Theorem, follows easily from the assumed regularity of $F$ and $\tilde{F}$. Especially for the finiteness of the RHS in \eqref{property}, we need to apply first the Cauchy-Schwartz inequality in $L^2\left(\left[0, \, T\right]\right)$ for the two Malliavin derivatives, and then, after an obvious cancellation, apply H\"older's inequality with the appropriate exponents to control the RHS by the product of two finite norms.

Taking now $\psi(y) = \mathbb{I}_{\left[a, \, b\right]}(y)$ for some $a, b \in \mathbb{R}$ with $a < b$ and $\phi(y) = \int_{-\infty}^{y}\psi(z)dz$, we can easily show that $\mathbb{P}$- almost surely we have $|\phi(F)| \leq b - a$, and also $D_{\cdot}\phi(F) = \psi(F)D_{\cdot}F$ by the comment after the proof of Proposition 1.2.3 on page 31 in \cite{Nualart} (since by our assumptions $D_{\cdot}F$ can never be identically zero, we can use Theorem 2.1.2 from page 86 in \cite{Nualart} to obtain absolute continuity). Thus, by the boundedness of $\psi(F)$ and our assumptions we have $\phi(F) \in L^{\infty}\left(\Omega\right) \cap \mathbb{D}^{1,2\tilde{\lambda} \vee \frac{q\tilde{r}}{q\tilde{r} - 1}}\left(\Omega\right)$, which is a subspace of $L^{\tilde{q}}\left(\Omega\right) \cap \mathbb{D}^{1, 2 \vee \frac{q\tilde{r}}{q\tilde{r} - 1}}\left(\Omega\right)$. Then, we can work as in the proof of Proposition 2.1.1 on page 78 in \cite{Nualart} to deduce that
\begin{eqnarray}\label{predensity}
\mathbb{E}\left[\psi(F)\right] = \mathbb{E}\left[\phi(F)\delta\left(\frac{D_{\cdot}F}{\left\Vert D_{\cdot}F \right\Vert_{L^2\left(\left[0, T\right]\right)}^2}\right)\right],
\end{eqnarray}
where now $\delta$ is the adjoint of the derivative operator 
\begin{equation}
D: L^{\tilde{q}}\left(\Omega\right) \cap \mathbb{D}^{1, 2 \vee \frac{q\tilde{r}}{q\tilde{r} - 1}}\left(\Omega \right) \longrightarrow L^{2\vee \frac{q\tilde{r}}{q\tilde{r} - 1}}\left(\Omega ; L^2\left(\left[0, \, T \right]\right)\right)
\end{equation}
(so an extension of the standard Skorokhod integral) the domain of which contains the process $\frac{D_{\cdot}F}{\left\Vert D_{\cdot}F \right\Vert_{L^2\left(\left[0, T\right]\right)}^2}$ as we have shown above. Since H\"older's inequality implies that
\begin{eqnarray}\label{CB}
\mathbb{E}\left[\left|\mathbb{I}_{F > x}\delta\left(\frac{D_{\cdot}F}{\left\Vert D_{\cdot}F \right\Vert_{L^2\left(\left[0, T\right]\right)}^2}\right)\right|\right] &\leq& \left \Vert \mathbb{I}_{F > x} \right \Vert_{L^{\tilde{q}}\left(\Omega\right)}\left \Vert \delta\left(\frac{D_{\cdot}F}{\left\Vert D_{\cdot}F \right\Vert_{L^2\left(\left[0, T\right]\right)}^2}\right) \right \Vert_{L^q\left(\Omega\right)} \nonumber \\
&\leq& \left \Vert \delta\left(\frac{D_{\cdot}F}{\left\Vert D_{\cdot}F \right\Vert_{L^2\left(\left[0, T\right]\right)}^2}\right) \right \Vert_{L^q\left(\Omega\right)}
\end{eqnarray}
which is finite, applying Fubini's Theorem on \eqref{predensity} we find that
\begin{eqnarray}
\mathbb{P}\left(a \leq F \leq b\right) = \int_{a}^{b}\mathbb{E}\left[\mathbb{I}_{F > x}\delta\left(\frac{D_{\cdot}F}{\left\Vert D_{\cdot}F \right\Vert_{L^2\left(\left[0, T\right]\right)}^2}\right)\right]dx.
\end{eqnarray}
This implies that $F$ has a density $f_F$ given by
\begin{eqnarray}
f_F(x) = \mathbb{E}\left[\mathbb{I}_{F > x}\delta\left(\frac{D_{\cdot}F}{\left\Vert D_{\cdot}F \right\Vert_{L^2\left(\left[0, T\right]\right)}^2}\right)\right] \nonumber
\end{eqnarray}
for all $x \in \mathbb{R}$, and by using the boundedness of the indicator function and the Dominated Convergence Theorem, we can show that this density is continuous. 

Finally, by recalling that the Skorokhod integral always has zero expectation (this also holds for its extension by a density argument), for $x \leq 0$ we have
\begin{eqnarray}\label{a1}
&& |x|^{\alpha}f_F(x) \nonumber \\
&&\qquad = \mathbb{E}\left[|x|^{\alpha}\mathbb{I}_{F > x}\delta\left(\frac{D_{\cdot}F}{\left\Vert D_{\cdot}F \right\Vert_{L^2\left(\left[0, T\right]\right)}^2}\right)\right] \nonumber \\
&&\qquad = -\mathbb{E}\left[|x|^{\alpha}\mathbb{I}_{F < x}\delta\left(\frac{D_{\cdot}F}{\left\Vert D_{\cdot}F \right\Vert_{L^2\left(\left[0, T\right]\right)}^2}\right)\right] \nonumber \\
&&\qquad = \mathbb{E}\left[|x|^{\alpha}\mathbb{I}_{F < -|x|}\left( -\delta\left(\frac{D_{\cdot}F}{\left\Vert D_{\cdot}F \right\Vert_{L^2\left(\left[0, T\right]\right)}^2}\right)\right)\right] \nonumber \\
&&\qquad \leq \mathbb{E}\left[|x|^{\alpha}\mathbb{I}_{|F| > |x|}\left|\delta\left(\frac{D_{\cdot}F}{\left\Vert D_{\cdot}F \right\Vert_{L^2\left(\left[0, T\right]\right)}^2}\right)\right|\right]
\end{eqnarray}
and for $x \geq 0$ obviously 
\begin{eqnarray}\label{a2}
&& |x|^{\alpha}f_F(x) \nonumber \\
&&\qquad = \mathbb{E}\left[|x|^{\alpha}\mathbb{I}_{F > x}\delta\left(\frac{D_{\cdot}F}{\left\Vert D_{\cdot}F \right\Vert_{L^2\left(\left[0, T\right]\right)}^2}\right)\right] \nonumber \\
&&\qquad = \mathbb{E}\left[|x|^{\alpha}\mathbb{I}_{F > |x|}\delta\left(\frac{D_{\cdot}F}{\left\Vert D_{\cdot}F \right\Vert_{L^2\left(\left[0, T\right]\right)}^2}\right)\right] \nonumber \\
&&\qquad \leq \mathbb{E}\left[|x|^{\alpha}\mathbb{I}_{|F| > |x|}\left|\delta\left(\frac{D_{\cdot}F}{\left\Vert D_{\cdot}F \right\Vert_{L^2\left(\left[0, T\right]\right)}^2}\right)\right|\right]
\end{eqnarray}
so by \eqref{a1}, \eqref{a2} and \eqref{densityest1} we obtain the estimate
\begin{eqnarray}
&& |x|^{\alpha}f_F(x) \nonumber \\
&&\qquad \leq \mathbb{E}\left[|F|^{\alpha}\mathbb{I}_{|F| > |x|}\left|\delta\left(\frac{D_{\cdot}F}{\left\Vert D_{\cdot}F \right\Vert_{L^2\left(\left[0, T\right]\right)}^2}\right)\right|\right] \nonumber \\
&&\qquad \leq \mathbb{E}\left[\left||F|^{\alpha}\delta\left(\frac{D_{\cdot}F}{\left\Vert D_{\cdot}F \right\Vert_{L^2\left(\left[0, T\right]\right)}^2}\right)\right|\right] \nonumber \\
&&\qquad \leq \left \Vert |F|^{\alpha}\delta\left(\frac{D_{\cdot}F}{\left\Vert D_{\cdot}F \right\Vert_{L^2\left(\left[0, T\right]\right)}^2}\right) \right \Vert_{L^q\left(\Omega\right)} \nonumber \\
&& \qquad \leq \left(C_{qr}+2\right)\mathbb{E}^{\frac{1}{qr}}\left[\left\Vert D_{\cdot, \cdot}^2F \right\Vert_{L^2\left(\left[0, T\right]^2\right)}^{qr}\right]\times\mathbb{E}^{\frac{1}{q\tilde{r}}}\left[\left|\frac{|F|^{\alpha}}{\left\Vert D_{\cdot}F \right\Vert_{L^2\left(\left[0, T\right]\right)}^2}\right|^{q\tilde{r}}\right] \nonumber \\
&& \qquad \qquad + C_{qr}\mathbb{E}^{\frac{1}{qr}}\left[\left\Vert D_{\cdot}F \right\Vert_{L^2\left(\left[0, T\right]\right)}^{qr}\right]\times\mathbb{E}^{\frac{1}{q\tilde{r}}}\left[\left|\frac{|F|^{\alpha}}{\left\Vert D_{\cdot}F \right\Vert_{L^2\left(\left[0, T\right]\right)}^2}\right|^{q\tilde{r}}\right] \nonumber
\end{eqnarray}
for all $x \in \mathbb{R}$. This completes the proof of the Lemma.
\end{proof}

\begin{proof}[\textbf{Proof of Lemma E2.2}]
We assume that the initial density $u_{0} = u_{0}(\cdot\,|\,\mathcal{G})$ is differentiable and that $w(\cdot)\left(u_{0}\right)_{x}$ is $L^2\left(\Omega \times \mathbb{R}^{+}\right)$-integrable. Then, by the theory developed in \cite{Krylov} we have that $u$ coincides with the unique solution to the SPDE 4.2 in a $w^2(\cdot)$ - weighted Sobolev space of higher regularity, and that \eqref{derivest} is also satisfied, with $M$ depending only on a compact interval $\mathcal{I} \subset \mathbb{R}^{+}$ which contains both $\displaystyle{\min_{0 \leq t \leq T}{\sigma_t}}$ and $\displaystyle{\max_{0 \leq t \leq T}{\sigma_t}}$. Note that even though the constants appearing in the Sobolev estimates obtained in \cite{Krylov} depend also on the modulus of spatial continuity of the coefficients of the SPDE, here this modulus of continuity is always zero since the coefficients do not depend on the spatial variable $x$. Next, we have
\begin{eqnarray}\label{formp}
\mathbb{E}\left[ \sup_{0 \leq t \leq T}\sup_{x \in \mathbb{R}^{+}}u^2(t,\,x) \right] \leq \mathbb{E}\left[ \sup_{0 \leq t \leq T}\sup_{x \in \left(0, \, 1\right)}u^2(t,\,x) \right] + \mathbb{E}\left[ \sup_{0 \leq t \leq T}\sup_{x \geq 1}u^2(t,\,x) \right], \nonumber \\
\end{eqnarray}
and we can use Morrey's inequality (see \cite{EVAN}) to control the second term in the RHS of the above by
\begin{eqnarray}
&&\mathbb{E}\left[ \sup_{0 \leq t \leq T}\int_{1}^{+\infty}u_{x}^2(t,\,x)dx \right] + \mathbb{E}\left[ \sup_{0 \leq t \leq T}\int_{1}^{+\infty}u^2(t,\,x)dx \right] \nonumber \\
&& \qquad \qquad \leq \mathbb{E}\left[ \sup_{0 \leq t \leq T}\int_{0}^{+\infty}w^2(x)u_{x}^2(t,\,x)dx \right] + \mathbb{E}\left[ \sup_{0 \leq t \leq T}\int_{0}^{+\infty}u^2(t,\,x)dx \right] \nonumber \\
&& \qquad \qquad \leq Me^{MT}\mathbb{E}\left[ \left \Vert w(\cdot)\left(u_0\right)_{x}(\cdot)\right \Vert_{L^{2}(\mathbb{R}^{+})}^{2}\right] + \left(Me^{MT} + 1\right)\mathbb{E}\left[ \left \Vert u_0(\cdot)\right \Vert_{L^{2}(\mathbb{R}^{+})}^{2}\right], \nonumber
\end{eqnarray}
where we have also used the identity (4.3) and estimate \eqref{derivest}. On the other hand, we can use Theorem 1 from \cite{DM13} to control the first term in the RHS of \eqref{formp} by 
\begin{equation}
\mathbb{E}\left[ \sup_{0 \leq t \leq T}u^2(t,\,1) \right] \leq \mathbb{E}\left[ \sup_{0 \leq t \leq T}\sup_{x \geq 1}u^2(t,\,x) \right] 
\end{equation}
which has already been controlled, and by the maximum of the initial density. Combining the above estimates we obtain \eqref{mp}.
\end{proof} 

\begin{proof}[\textbf{Proof of Theorem E1.1}]
By Lemmas 3.2, 3.3 and 3.4, we have that $\sigma_{t}$ satisfies the assumptions of Lemma~\ref{lem:3.1} for any $q, r > 1$ with $qr < \frac{4k\theta}{3\xi^2}$, any $\alpha \geq 0$, and any $\lambda \geq q\tilde{r}$, under the conditional probability measure $\mathbb{P}(\cdot\,|\,B_{\cdot}^{0},\,\mathcal{G})$, since we can show that $q, r$ can be chosen such that $\frac{\sigma_{t}^{\alpha}}{\left\Vert D_{\cdot}\sigma_t \right\Vert_{L^2\left(\left[0, T\right]\right)}^2} \in L_{B_{\cdot}^{0},\,\mathcal{G}}^{q\tilde{r}}\left(\Omega\right)$, $\mathbb{P}$-almost surely. To see the last, we use the Cauchy-Schwartz inequality and we recall Theorem~3.1 from \cite{HK06} to obtain
\begin{eqnarray}\label{cirfin}
&&\mathbb{E}\left[\frac{\sigma_{t}^{q\tilde{r}\alpha}}{\left\Vert D_{\cdot}\sigma_t \right\Vert_{L^2\left(\left[0, T\right]\right)}^{2q\tilde{r}}}\right] \nonumber \\
&& \qquad \qquad\leq \frac{\xi\sqrt{1-\rho^{2}_{2}}}{t^{2q\tilde{r}}}\mathbb{E}\left[\sigma_{t}^{q\tilde{r}(\alpha - 1)}\left(\int_{0}^{t}e^{2\int_{t'}^{t}\left[\left(\frac{k\theta}{2}-\frac{\xi^{2}}{8}\right)\frac{1}{\sigma_{s}}+\frac{k}{2}\right]ds}\right)^{q\tilde{r}}dt'\right] \nonumber \\
&& \qquad \qquad\leq\frac{\xi\sqrt{1-\rho^{2}_{2}}}{t^{q\tilde{r}}}\mathbb{E}\left[\sigma_{t}^{q\tilde{r}(\alpha - 1)}e^{2q\tilde{r}\int_{0}^{t}\left[\left(\frac{k\theta}{2}-\frac{\xi^{2}}{8}\right)\frac{1}{\sigma_{s}}+\frac{k}{2}\right]ds}\right] \nonumber \\
&& \qquad \qquad= \frac{\xi\sqrt{1-\rho^{2}_{2}}}{t^{q\tilde{r}}}\mathbb{E}\left[\mathbb{E}\left[\sigma_{t}^{q\tilde{r}(\alpha - 1)}e^{2q\tilde{r}\int_{0}^{t}\left[\left(\frac{k\theta}{2}-\frac{\xi^{2}}{8}\right)\frac{1}{\sigma_{s}}+\frac{k}{2}\right]ds}\, \Big| \, \sigma_0\right]\right] \nonumber \\
&& \qquad \qquad\leq \frac{\tilde{C}}{t^{q\tilde{r}}}\gamma_t^{v - q\tilde{r}(\alpha - 1)}\mathbb{E}\left[\sigma_{0}^{v}H\left(-\gamma_t\sigma_0e^{-kt}\right)\right]
\end{eqnarray}
for some $c, \tilde{C} > 0$, provided that $2q\tilde{r}\left(\frac{k\theta}{2}-\frac{\xi^{2}}{8}\right) < \frac{\xi^2}{8}\left(\frac{2k\theta}{\xi^2} - 1\right)^2$ and $q\tilde{r}(\alpha - 1) > -\frac{2k\theta}{\xi^2} - v$, where
\begin{eqnarray}
v = \frac{1}{2}\left(-\left(\frac{2k\theta}{\xi^2} - 1\right) + \sqrt{\left(\frac{2k\theta}{\xi^2} - 1\right)^2 - 2q\tilde{r}\left(\frac{4k\theta}{\xi^2}-1\right)}\right),
\end{eqnarray}
$\gamma_{t}=\frac{2k}{\xi^{2}}\left(1-e^{-kt}\right)^{-1}>\frac{2k}{\xi^{2}}$
for all $t\geq0$, and $H$ is a hypergeometric function for which
we have the asymptotic estimate of page 17 in \cite{HK06}.
To have $2q\tilde{r}\left(\frac{k\theta}{2}-\frac{\xi^{2}}{8}\right) < \frac{\xi^2}{8}\left(\frac{2k\theta}{\xi^2} - 1\right)^2$, $q\tilde{r}(\alpha - 1) > -\frac{2k\theta}{\xi^2} - v$ and $qr < \frac{4k\theta}{3\xi^2}$ for sufficiently small $q > 1$, it suffices to obtain all these strict inequalities for $q = 1$. Since $r < \frac{4k\theta}{3\xi^2}$ is equivalent to $\tilde{r} > \frac{4x}{4x - 3}$ for $x = \frac{k\theta}{\xi^2}$ with $\frac{4x}{4x - 3} > 1$, we can have this inequality along with $2\tilde{r}\left(\frac{k\theta}{2}-\frac{\xi^{2}}{8}\right) < \frac{\xi^2}{8}\left(\frac{2k\theta}{\xi^2} - 1\right)^2 \Leftrightarrow \tilde{r} < \frac{(2x - 1)^2}{2(4x - 1)}$ for some $\tilde{r} > 1$ if and only if $\frac{4x}{4x - 3} <  \frac{(2x - 1)^2}{2(4x - 1)}$. The last inequality is satisfied since it is equivalent to $16x^3 - 60x^2 + 24x - 3 > 0$ and we have $x = \frac{k\theta}{\xi^{2}} > x^{*}$. Then, $q\tilde{r}(\alpha - 1) > -\frac{2k\theta}{\xi^2} - v$ can be obtained for $p = 1$, for any $\alpha \geq 0$ and $\tilde{r}$ sufficiently close to its upper bound $\frac{(2x - 1)^2}{2(4x - 1)}$, provided that it holds for $q = 1, \alpha = 0$ and $\tilde{r} = \frac{(2x - 1)^2}{2(4x - 1)}$, i.e when
\begin{eqnarray}\label{weightest2}
-\frac{(2x - 1)^2}{2(4x - 1)} > -2x + \frac{1}{2}\left(2x - 1\right)
\end{eqnarray}
which is also satisfied when $x = \frac{k\theta}{\xi^{2}} > x^{*}$ (since $x^{*} > 1$). Next, since $\sigma_0$ is bounded away from zero, the argument of $H$ in \eqref{cirfin} is bounded from below by some $m > 0$, and then the estimate of page 17 in \cite{HK06} gives $H(-z)\leq K|z|^{-v + q\tilde{r}(\alpha - 1)}$ for all $z \geq m$, for some $K > 0$. Therefore, from \eqref{cirfin} we obtain 
\begin{eqnarray}\label{cirfin2}
\mathbb{E}\left[\frac{\sigma_{t}^{q\tilde{r}\alpha}}{\left\Vert D_{\cdot}\sigma_t \right\Vert_{L^2\left(\left[0, T\right]\right)}^{2q\tilde{r}}}\right] \leq \frac{\tilde{C}}{t^{q\tilde{r}}}e^{cT}\mathbb{E}\left[\sigma_{0}^{q\tilde{r}(\alpha - 1)}\right]
\end{eqnarray}
for some $c > 0$, with the RHS of the above being finite, and this implies also 
\begin{eqnarray}
\mathbb{E}\left[\frac{\sigma_{t}^{q\tilde{r}\alpha}}{\left\Vert D_{\cdot}\sigma_t \right\Vert_{L^2\left(\left[0, T\right]\right)}^{2q\tilde{r}}}\,|\,B_{\cdot}^{0},\,\mathcal{G}\right] < \infty
\end{eqnarray}
$\mathbb{P}$- almost surely. The last means that the assumptions of Lemma~\ref{lem:3.1} are indeed satisfied. 

From the above we deduce that under $\mathbb{P}(\cdot\,|\,B_{\cdot}^{0},\,\mathcal{G})$, $\sigma_t$ has a density $p_{t}(y\,|\,B_{\cdot}^{0},\,\mathcal{G})$ which is supported in $\left[0, \, +\infty\right)$ (since this CIR process does not hit zero) and which satisfies
\begin{eqnarray}\label{3.2est}
&& \sup_{y \in \mathbb{R}^{+}}y^{\alpha}p_{t}(y\,|\,B_{\cdot}^{0},\,\mathcal{G}) \nonumber \\
&& \qquad \leq \left(C+2\right)\mathbb{E}^{\frac{1}{qr}}\left[\left\Vert D_{\cdot, \cdot}^2\sigma_t \right\Vert_{L^2\left(\left[0, T\right]^2\right)}^{qr}\,|\,B_{\cdot}^{0},\,\mathcal{G}\right] \nonumber \\
&& \qquad \qquad \qquad \times\mathbb{E}^{\frac{1}{q\tilde{r}}}\left[\left|\frac{\sigma_t^{\alpha}}{\left\Vert D_{\cdot}\sigma_t \right\Vert_{L^2\left(\left[0, T\right]\right)}^2}\right|^{q\tilde{r}}\,|\,B_{\cdot}^{0},\,\mathcal{G}\right] \nonumber \\
&& \qquad \qquad + C\mathbb{E}^{\frac{1}{qr}}\left[\left\Vert D_{\cdot}\sigma_t \right\Vert_{L^2\left(\left[0, T\right]\right)}^{qr}\,|\,B_{\cdot}^{0},\,\mathcal{G}\right]\times\mathbb{E}^{\frac{1}{q\tilde{r}}}\left[\left|\frac{\sigma_t^{\alpha}}{\left\Vert D_{\cdot}\sigma_t \right\Vert_{L^2\left(\left[0, T\right]\right)}^2}\right|^{q\tilde{r}} \,|\,B_{\cdot}^{0},\,\mathcal{G}\right] \nonumber \\
&& \qquad = \left(C+2\right)\mathbb{E}^{\frac{1}{qr}}\left[\left\Vert D_{\cdot, \cdot}^2\sigma_t \right\Vert_{L^2\left(\left[0, t\right]^2\right)}^{qr}\,|\,B_{\cdot}^{0},\,\mathcal{G}\right] \nonumber \\
&& \qquad \qquad \qquad \times\mathbb{E}^{\frac{1}{q\tilde{r}}}\left[\left|\frac{\sigma_t^{\alpha}}{\left\Vert D_{\cdot}\sigma_t \right\Vert_{L^2\left(\left[0, t\right]\right)}^2}\right|^{q\tilde{r}}\,|\,B_{\cdot}^{0},\,\mathcal{G}\right] \nonumber \\
&& \qquad \qquad + C\mathbb{E}^{\frac{1}{qr}}\left[\left\Vert D_{\cdot}\sigma_t \right\Vert_{L^2\left(\left[0, t\right]\right)}^{qr}\,|\,B_{\cdot}^{0},\,\mathcal{G}\right]\times\mathbb{E}^{\frac{1}{q\tilde{r}}}\left[\left|\frac{\sigma_t^{\alpha}}{\left\Vert D_{\cdot}\sigma_t \right\Vert_{L^2\left(\left[0, t\right]\right)}^2}\right|^{q\tilde{r}}\,|\,B_{\cdot}^{0},\,\mathcal{G}\right], \nonumber 
\end{eqnarray}
so raising to the power $q$, taking expectations and using Holder's inequality we obtain
\begin{eqnarray}\label{3.2est2}
&& \mathbb{E}\left[\left(\sup_{y \in \mathbb{R}^{+}}y^{\alpha}p_{t}(y\,|\,B_{\cdot}^{0},\,\mathcal{G})\right)^q\right] \nonumber \\
&& \qquad \leq \tilde{C}\mathbb{E}^{\frac{1}{r}}\left[\left\Vert D_{\cdot, \cdot}^2\sigma_t \right\Vert_{L^2\left(\left[0, t\right]^2\right)}^{qr}\right]\times\mathbb{E}^{\frac{1}{\tilde{r}}}\left[\left|\frac{\sigma_t^{\alpha}}{\left\Vert D_{\cdot}\sigma_t \right\Vert_{L^2\left(\left[0, t\right]\right)}^2}\right|^{q\tilde{r}}\right] \nonumber \\
&& \qquad \qquad + \tilde{C}\mathbb{E}^{\frac{1}{r}}\left[\left\Vert D_{\cdot}\sigma_t \right\Vert_{L^2\left(\left[0, t\right]\right)}^{qr}\right]\times\mathbb{E}^{\frac{1}{\tilde{r}}}\left[\left|\frac{\sigma_t^{\alpha}}{\left\Vert D_{\cdot}\sigma_t \right\Vert_{L^2\left(\left[0, t\right]\right)}^2}\right|^{q\tilde{r}}\right]
\end{eqnarray}
for some $\tilde{C} > 0$. Next, by Lemmas 3.2 and 3.3 we have
\begin{eqnarray}\label{subst1}
\mathbb{E}^{\frac{1}{qr}}\left[\left\Vert D_{\cdot}\sigma_t \right\Vert_{L^2\left(\left[0, t\right]\right)}^{qr}\right] &=& \mathbb{E}^{\frac{1}{qr}}\left[\left(\int_{0}^{t}\xi^2\left(1-\rho^{2}_{2}\right)e^{-2\int_{t'}^{t}\left[\left(\frac{k\theta}{2}-\frac{\xi^{2}}{8}\right)\frac{1}{\sigma_{s}}+\frac{k}{2}\right]ds}\sigma_{t}dt'\right)^{\frac{qr}{2}}\right] \nonumber \\
&\leq& \xi\sqrt{1-\rho^{2}_{2}}\mathbb{E}^{\frac{1}{qr}}\left[\left(\int_{0}^{t}\sup_{s \leq T}\sigma_{s}dt'\right)^{\frac{qr}{2}}\right] \nonumber \\
&=& \xi\sqrt{1-\rho^{2}_{2}}\sqrt{t}\mathbb{E}^{\frac{1}{qr}}\left[\left(\sup_{s \leq T}\sigma_{s}\right)^{\frac{qr}{2}}\right] \nonumber \\
&=& C'\sqrt{t}
\end{eqnarray}
for some $C' > 0$, while by the estimate of Lemma 3.4 for $q' = qr$ we have
\begin{eqnarray}\label{subst2}
\mathbb{E}^{\frac{1}{qr}}\left[\left\Vert D_{\cdot, \cdot}^2\sigma_t \right\Vert_{L^2\left(\left[0, t\right]\right)}^{qr}\right] \leq \mathbb{E}^{\frac{1}{qr}}\left[t^{qr}\sup_{0\leq t',t''\leq t\leq T}|D_{t',t''}^{2}\sigma_{t}|^{qr}\right] = C''t
\end{eqnarray}
for some $C'' > 0$, since $qr < \frac{4k\theta}{3\xi^2}$. Moreover, for our choice of $q$ and $r$, by \eqref{cirfin2} we have
\begin{eqnarray}\label{subst3}
\mathbb{E}^{\frac{1}{q\tilde{r}}}\left[\frac{\sigma_{t}^{q\tilde{r}\alpha}}{\left\Vert D_{\cdot}\sigma_t \right\Vert_{L^2\left(\left[0, T\right]\right)}^{2q\tilde{r}}}\right] \leq C^{(3)}\frac{1}{t}
\end{eqnarray}
for some $C^{(3)} > 0$. Substituting now \eqref{subst1}, \eqref{subst2} and \eqref{subst3} in \eqref{3.2est2}, we obtain 
\begin{eqnarray}
&& \mathbb{E}\left[\left(\sup_{y \in \mathbb{R}^{+}}y^{\alpha}p_{t}(y\,|\,B_{\cdot}^{0},\,\mathcal{G}^0)\right)^q\right] \leq C^{(4)} + C^{(5)}\frac{1}{\sqrt{t^q}}
\end{eqnarray}
for some $C^{(4)}, C^{(5)} > 0$. Since the RHS of the last is integrable in $t$ for $t \in \left[0, \, T\right]$ (since we can take $q < 2$), the desired result follows.
\end{proof}

\begin{proof}[\textbf{Proof of Theorem E1.2}]
Let $f$ be a smooth function, compactly supported
in $\mathbb{R}^{2}$, such that $f$ vanishes on the $y$ - axis. Theorem~\ref{thm:3.2} applied on the $\left(W_{\cdot}^1, W_{\cdot}^0\right)$ - driven CIR process $\left\{\sigma_{t}^1 : t \geq 0\right\}$ implies that the last possesses a density $p_{t}\left(\cdot | B_{\cdot}^0, \mathcal{G}\right)$ for each $t \geq 0$, for which we have
\begin{eqnarray}
v_{t,C_{1}}\left(f \right)&=& \mathbb{E}\left[f\left(X_{t}^{1}, \, \sigma_{t}^{1}\right)\mathbb{I}_{\{T_{1}\geq t\}}\,|W_{\cdot}^{0},B_{\cdot}^{0},C_{1},\mathcal{G}\right] \nonumber \\
&=&\mathbb{E}\left[\mathbb{E}\left[f\left(X_{t}^{1}, \, \sigma_{t}^{1}\right)\mathbb{I}_{\{T_{1}\geq t\}}\,|W_{\cdot}^{0},\sigma_{t}^{1},B_{\cdot}^{0},C_{1},\mathcal{G}\right]\,|W_{\cdot}^{0},B_{\cdot}^{0},C_{1},\mathcal{G}\right]
\nonumber \\
&=&\int_{\mathbb{R}}\mathbb{E}\left[f\left(X_{t}^{1}, \, y\right)\mathbb{I}_{\{T_{1}\geq t\}}\,|W_{\cdot}^{0},\sigma_{t}^{1}=y,B_{\cdot}^{0},C_{1},\mathcal{G}\right]p_{t}\left(y|B_{\cdot}^{0},\mathcal{G}\right)dy. \label{eq:4.9}
\end{eqnarray}
for any $t \geq 0$. Next, we compute
\begin{eqnarray}
& & \mathbb{E}\left[f\left(X_{t}^{1}, \, y\right)\mathbb{I}_{\{T_{1}\geq t\}}|W_{\cdot}^{0},
\sigma_{t}^{1}=y,B_{\cdot}^{0},C_{1},\mathcal{G}\right] \nonumber \\
& & \qquad =\mathbb{E}\left[\mathbb{E}\left[f\left(X_{t}^{1}, \, y\right)\mathbb{I}_{\{T_{1}\geq t\}}|
W_{\cdot}^{0},\sigma_{.},C_{1},\mathcal{G}\right]|W_{\cdot}^{0},\sigma_{t}^{1}=y,B_{\cdot}^{0},C_{1},
\mathcal{G}\right] \nonumber \\
& & \qquad =\mathbb{E}\left[\int_{\mathbb{R}^+}f(x, \, y)u\left(t,x,W_{\cdot}^{0},\mathcal{G},C_{1}, 
h\left(\sigma_{.}^1\right)\right)dx|W_{\cdot}^{0},\sigma_{t}^{1}=y,B_{\cdot}^{0},C_{1},\mathcal{G}\right] 
\nonumber \\
& & \qquad =\int_{\mathbb{R}^+}f(x, \, y)\mathbb{E}\left[u\left(t,x,W_{\cdot}^{0},\mathcal{G},C_{1}, 
h\left(\sigma_{.}^1\right)\right) |W_{\cdot}^{0},\sigma_{t}^{1}=y,B_{\cdot}^{0},C_{1},\mathcal{G}\right]dx, \nonumber \\ \label{eq:4.10}
\end{eqnarray}
where $u\left(t,x,W_{\cdot}^0,C_{1},\mathcal{G},h\left(\sigma_{.}^1\right)\right)$
is the $L^{2}\left(\Omega\times[0,\,T];\,H_{0}^{1}\left(\mathbb{R}^{+}\right)\right)$
density given by Theorem~4.1 when the coefficient vector $C_{1}$
is given and the volatility path is $h\left(\sigma_{.}^1\right)$. By \eqref{eq:4.9} and \eqref{eq:4.10} 
we have that the desired density exists and is given by
\begin{eqnarray}
&&u_{C_{1}}\left(t,x,y,W_{\cdot}^{0},B_{\cdot}^{0},\mathcal{G}\right) \nonumber \\
&& \qquad \qquad = p_{t}\left(y|B_{\cdot}^{0},\mathcal{G}\right)
\mathbb{E}\left[u\left(t,x,W_{\cdot}^{0},\mathcal{G},C_{1},h\left(\sigma_{.}^1\right)\right)\,|W_{\cdot}^{0},\sigma_{t}^{1}=
y,B_{\cdot}^{0},C_{1},\mathcal{G}\right] \nonumber \\
\end{eqnarray}
which is supported in $\mathbb{R}^{+}\times\mathbb{R}^{+}$. Using the Cauchy-Schwartz inequality, the law of total expectation, Fubini's Theorem, and the identity (4.3) we obtain for any $\alpha \geq 0$
\begin{eqnarray*}
& & \int_{\mathbb{R}^{+}}\int_{\mathbb{R}^{+}}y^{a}\left(u_{C_{1}}\left(t,\,x,\,y,\,W_{\cdot}^{0},\,B_{\cdot}^{0},\,\mathcal{G}\right)\right)^{2}dydx \\
& & \qquad = \int_{\mathbb{R}^{+}}\int_{\mathbb{R}^{+}}y^{\alpha}p_{t}^2\left(y|B_{\cdot}^{0},\mathcal{G}\right)
\mathbb{E}^2\left[u\left(t,x,W_{\cdot}^{0},\mathcal{G},C_{1},h\left(\sigma_{.}^1\right)\right)\,|
W_{\cdot}^{0},B_{\cdot}^{0},\sigma_{t}^{1}=
y,C_{1},\mathcal{G}\right]dydx \\
& & \qquad \leq M_{B_{\cdot}^{0},\, \mathcal{G}}^{\alpha}(t)\int_{\mathbb{R}^{+}}\int_{\mathbb{R}^{+}}p_{t}\left(y|B_{\cdot}^{0},\mathcal{G}\right) \nonumber \\
& & \qquad \qquad \times \mathbb{E}\left[u^2\left(t,x,W_{\cdot}^{0},\mathcal{G},C_{1},h\left(\sigma_{.}^1\right)\right)\,|
W_{\cdot}^{0},B_{\cdot}^{0},\sigma_{t}^{1}=
y,C_{1},\mathcal{G}\right]dydx \\
& & \qquad = M_{B_{\cdot}^{0},\, \mathcal{G}}^{\alpha}(t)\int_{\mathbb{R}^{+}}
\mathbb{E}\left[u^{2}\left(t,x,W_{\cdot}^{0},\mathcal{G},C_{1},h\left(\sigma_{.}^1\right)\right)\,|
W_{\cdot}^{0},B_{\cdot}^{0},C_{1},\mathcal{G}\right]dx \\
& & \qquad = M_{B_{\cdot}^{0},\, \mathcal{G}}^{\alpha}(t)
\mathbb{E}\left[\int_{\mathbb{R}^{+}}u^{2}\left(t,x,W_{\cdot}^{0},\mathcal{G},C_{1},h\left(\sigma_{.}^1\right)\right)dx\,|
W_{\cdot}^{0},B_{\cdot}^{0},C_{1},\mathcal{G}\right] \\
& & \qquad \leq M_{B_{\cdot}^{0},\, \mathcal{G}}^{\alpha}(t)\mathbb{E}\left[ \left \Vert u_0(\cdot)\right \Vert_{L^{2}(\mathbb{R}^{+})}^{2}\,|\,\mathcal{G}\right],
\end{eqnarray*}
where we have $\displaystyle{M_{B_{\cdot}^{0},\, \mathcal{G}}^{\alpha}(\cdot) = \sup_{y\geq0}\left(y^{\alpha}p_{\cdot}(y\,|\,B_{\cdot}^{0},\,\mathcal{G})\right) \in L^{q}\left(\Omega \times \left[0, \, T\right]\right)}$
for all small enough $q > 1$ (given $C_1$, by Theorem~\ref{thm:3.2}). Denoting by $\mathbb{E}_{C_1}$ the expectation given $C_1$ and taking $q' > 1$ such that $\frac{1}{q} + \frac{1}{q'} = 1$, by the above and by Holder's inequality we get
\begin{eqnarray*}
& & \mathbb{E}_{C_1}\left[\int_{0}^{T}\int_{\mathbb{R}^{+}}\int_{\mathbb{R}^{+}}y^{a}\left(u_{C_{1}}\left(t,\,x,\,y,\,W_{\cdot}^{0},\,B_{\cdot}^{0},\,\mathcal{G}\right)\right)^{2}dydxdt\right] \\
& & \qquad \leq \mathbb{E}_{C_1}\left[\int_{0}^{T}M_{B_{\cdot}^{0},\, \mathcal{G}}^{\alpha}(t)\mathbb{E}\left[ \left \Vert u_0(\cdot)\right \Vert_{L^{2}(\mathbb{R}^{+})}^{2}\,|\,\mathcal{G}\right]dt\right] \\
& & \qquad \leq \mathbb{E}_{C_1}^{\frac{1}{q}}\left[\left(\int_{0}^{T}M_{B_{\cdot}^{0},\, \mathcal{G}}^{\alpha}(t)dt\right)^q\right] \mathbb{E}_{C_1}^{\frac{1}{q'}}\left[\mathbb{E}^{q'}\left[ \left \Vert u_0(\cdot)\right \Vert_{L^{2}(\mathbb{R}^{+})}^{2}\,|\,\mathcal{G}\right]\right] \\
& & \qquad \leq T^{\frac{1}{q'}}\mathbb{E}_{C_1}^{\frac{1}{q}}\left[\int_{0}^{T}\left(M_{B_{\cdot}^{0},\, \mathcal{G}}^{\alpha}(t)\right)^q dt\right] \mathbb{E}^{\frac{1}{q'}}\left[\mathbb{E}^{q'}\left[ \left \Vert u_0(\cdot)\right \Vert_{L^{2}(\mathbb{R}^{+})}^{2}\,|\,\mathcal{G}\right]\right] \\
& & \qquad < \infty,
\end{eqnarray*}
which shows that the density belongs to the space $L_{\alpha}$ for any $\alpha \geq 0$. Moreover, repeating the above computations but for the derivative multiplied by $w(x)$, we find that
\begin{eqnarray*}
& & \int_{\mathbb{R}^{+}}\int_{\mathbb{R}^{+}}w^2(x)y^{a}\left(\frac{\partial u_{C_{1}}}{\partial x}\left(t,\,x,\,y,\,W_{\cdot}^{0},\,B_{\cdot}^{0},\,\mathcal{G}\right)\right)^{2}dydx \\
& & \qquad \leq M_{B_{\cdot}^{0},\, \mathcal{G}}^{\alpha}(t)\mathbb{E}\left[\int_{\mathbb{R}^{+}}w^2(x)u_{x}^{2}\left(t,\,x,\,W_{\cdot}^{0},\,\mathcal{G},\,C_{1},\,h\left(\sigma_{.}^1\right)\right)dx\,|\,W_{\cdot}^{0},\,B_{\cdot}^{0},\,C_{1},\,\mathcal{G}\right],
\end{eqnarray*}
so when $\rho_3 := \int_{0}^{1}dW^0_tdB^0_t = 0$, writing $\mathbb{E}_{C_{1}, \mathcal{G}}^{B_{\cdot}^{0}}$ for the expectation given $C_{1}$, $\mathcal{G}$ and $B_{\cdot}^{0}$ and using Lemma~\ref{lem:4dest}, we obtain
\begin{eqnarray*}
& & \mathbb{E}_{C_{1}, \mathcal{G}}^{B_{\cdot}^{0}}\left[\int_{\mathbb{R}^{+}}\int_{\mathbb{R}^{+}}w^2(x)y^{a}\left(\frac{\partial u_{C_{1}}}{\partial x}\left(t,\,x,\,y,\,W_{\cdot}^{0},\,B_{\cdot}^{0},\,\mathcal{G}\right)\right)^{2}dydx \right] \\
& & \qquad \leq M_{B_{\cdot}^{0},\, \mathcal{G}}^{\alpha}(t)\mathbb{E}\left[\int_{\mathbb{R}^{+}}w^2(x)u_{x}^{2}\left(t,\,x,\,W_{\cdot}^{0},\,\mathcal{G},\,C_{1},\,h\left(\sigma_{.}^1\right)\right)dx\,|\,B_{\cdot}^{0},\,C_{1},\,\mathcal{G}\right] \\
& & \qquad \leq M_{B_{\cdot}^{0},\, \mathcal{G}}^{\alpha}(t)\mathbb{E}\left[\sup_{0 \leq s \leq T}\int_{\mathbb{R}^{+}}w^2(x)u_{x}^{2}\left(s,\,x,\,W_{\cdot}^{0},\,\mathcal{G},\,C_{1},\,h\left(\sigma_{.}^1\right)\right)dx\,|\,B_{\cdot}^{0},\,C_{1},\,\mathcal{G}\right] \\
& & \qquad \leq Me^{MT}M_{B_{\cdot}^{0},\, \mathcal{G}}^{\alpha}(t)\left(\mathbb{E}\left[ \left \Vert w(\cdot)\left(u_0\right)_{x}(\cdot)\right \Vert_{L^{2}(\mathbb{R}^{+})}^{2}\,|\,\mathcal{G}\right] + \mathbb{E}\left[\left \Vert u_0(\cdot)\right \Vert_{L^{2}(\mathbb{R}^{+})}^{2}\,|\,\mathcal{G}\right]\right).
\end{eqnarray*}
The last implies that 
\begin{eqnarray*}
& & \mathbb{E}_{C_{1}}\left[\int_{0}^{T}\int_{\mathbb{R}^{+}}\int_{\mathbb{R}^{+}}w^2(x)y^{a}\left(\frac{\partial u_{C_{1}}}{\partial x}\left(t,\,x,\,y,\,W_{\cdot}^{0},\,B_{\cdot}^{0},\,\mathcal{G}\right)\right)^{2}dydxdt \right] \\
& & \qquad \leq Me^{MT}\mathbb{E}_{C_{1}}\left[\int_{0}^{T}M_{B_{\cdot}^{0},\, \mathcal{G}}^{\alpha}(t)dt\mathbb{E}\left[ \left \Vert w(\cdot)\left(u_0\right)_{x}(\cdot)\right \Vert_{L^{2}(\mathbb{R}^{+})}^{2}\,|\,\mathcal{G}\right]\right] \\
& & \qquad \qquad + Me^{MT}\mathbb{E}_{C_{1}}\left[\int_{0}^{T}M_{B_{\cdot}^{0},\, \mathcal{G}}^{\alpha}(t)dt\mathbb{E}\left[\left \Vert u_0(\cdot)\right \Vert_{L^{2}(\mathbb{R}^{+})}^{2}\,|\,\mathcal{G}\right]\right] \\
& & \qquad \leq MT^{\frac{1}{q'}}e^{MT}\mathbb{E}_{C_1}^{\frac{1}{q}}\left[\int_{0}^{T}\left(M_{B_{\cdot}^{0},\, \mathcal{G}}^{\alpha}(t)\right)^q dt\right] \nonumber \\
& & \qquad \qquad \qquad \times \mathbb{E}_{C_1}^{\frac{1}{q'}}\left[\mathbb{E}^{q'}\left[ \left \Vert w(\cdot)\left(u_0\right)_{x}(\cdot)\right \Vert_{L^{2}(\mathbb{R}^{+})}^{2}\,|\,\mathcal{G}\right]\right] \\
& & \qquad \qquad + MT^{\frac{1}{q'}}e^{MT}\mathbb{E}_{C_1}^{\frac{1}{q}}\left[\int_{0}^{T}\left(M_{B_{\cdot}^{0},\, \mathcal{G}}^{\alpha}(t)\right)^q dt\right] \mathbb{E}_{C_1}^{\frac{1}{q'}}\left[\mathbb{E}^{q'}\left[ \left \Vert u_0(\cdot)\right \Vert_{L^{2}(\mathbb{R}^{+})}^{2}\,|\,\mathcal{G}\right]\right] \\
& & \qquad < \infty
\end{eqnarray*}
which gives the weighted integrability of the derivative. To obtain the boundary condition when $\rho_3 = 0$, we work as follows
\begin{eqnarray*}
& & \mathbb{E}_{C_1}\left[\int_{\mathbb{R}^{+}}\int_{\mathbb{R}^{+}}y^{\alpha}\left(u_{C_{1}}\left(t,\,x,\,y,\,W_{\cdot}^{0},\,B_{\cdot}^{0},\,\mathcal{G}\right)\right)^{2}dydt\right] \\
& & \qquad \leq \mathbb{E}_{C_1}\Bigg[\int_{\mathbb{R}^{+}}\int_{\mathbb{R}^{+}}M_{B_{\cdot}^{0},\, \mathcal{G}}^{\alpha}(t) p_{t}\left(y\,|\,B_{\cdot}^{0},\,\mathcal{G}\right) \\
& & \qquad \qquad \times \mathbb{E}^{2}\left[u\left(t,\,x,\,W_{\cdot}^{0},\,\mathcal{G},\,C_{1},\,h\left(\sigma_{.}^1\right)\right)\,|\,W_{\cdot}^{0},\,\sigma_{t}^{1}=y,\,B_{\cdot}^{0},\,C_{1},\,\mathcal{G}\right]dydt\Bigg] \\
& & \qquad \leq \mathbb{E}_{C_1}\Bigg[\int_{\mathbb{R}^{+}}M_{B_{\cdot}^{0},\, \mathcal{G}}^{\alpha}(t)\int_{\mathbb{R}^{+}}p_{t}\left(y\,|\,B_{\cdot}^{0},\,\mathcal{G}\right) \nonumber \\
& & \qquad \qquad \times \mathbb{E}\left[u^{2}\left(t,\,x,\,W_{\cdot}^{0},\,\mathcal{G},\,C_{1},\,h\left(\sigma_{.}^1\right)\right)\,|\,W_{\cdot}^{0},\,\sigma_{t}^{1}=y,\,B_{\cdot}^{0},\,C_{1},\,\mathcal{G}\right]dydt\Bigg] \\
& & \qquad = \mathbb{E}_{C_1}\left[\int_{\mathbb{R}^{+}}M_{B_{\cdot}^{0},\, \mathcal{G}}^{\alpha}(t)\mathbb{E}\left[u^{2}\left(t,\,x,\,W_{\cdot}^{0},\,\mathcal{G},\,C_{1},\,h\left(\sigma_{.}^1\right)\right)\,|\,W_{\cdot}^{0},\,B_{\cdot}^{0},\,C_{1},\,\mathcal{G}\right]dt\right] \\
& & \qquad = \mathbb{E}_{C_1}\left[\int_{\mathbb{R}^{+}}M_{B_{\cdot}^{0},\, \mathcal{G}}^{\alpha}(t)\mathbb{E}_{C_1}\left[u^{2}\left(t,\,x,\,W_{\cdot}^{0},\,\mathcal{G},\,C_{1},\,h\left(\sigma_{.}^1\right)\right)\,|\,B_{\cdot}^{0},\,C_{1},\,\mathcal{G}\right]dt\right],
\end{eqnarray*}
where we can use the maximum principle given in Lemma~\ref{lem:4dest}, the integrability of $M_{B_{\cdot}^{0},\, \mathcal{G}}^{\alpha}(\cdot)$ and the Dominated Convergence Theorem, to show that the RHS of the last tends to zero as $x \longrightarrow 0^{+}$. This completes the proof of the Theorem.
\end{proof}

\begin{proof}[\textbf{Proof of Lemma E2.3}]
The finiteness of all the terms in the identity we are proving is a consequence of Lemma~5.3 and the assumed weighted integrability of $u$ and $u_{x}$. Multiplying equation (5.6) by $w^2(x)\left(y^{\delta}\right)^{+}$, applying Ito's formula for the $L^{2}(\mathbb{R}^{+})$ norm (Theorem~3.1 from \cite{KR81} for the triple $H_{0}^{1}\subset L^{2}\subset H^{-1}$, with $\Lambda(u) = w(\cdot)u$), and then integrating in
$y$ over $\mathbb{R}^{+}$, we obtain the equality
\begin{eqnarray} 
& & \left\Vert I_{\epsilon,1}(t,\,\cdot)\right\Vert _{\tilde{L}_{\delta, w}^2}^{2}=\left\Vert \int_{\mathcal{D}}U_{0}(\cdot,z)
\phi_{\epsilon}(z,\cdot)dz\right\Vert _{\tilde{L}_{\delta, w}^2}^{2} \nonumber \\
& & \qquad -2r_{1}\int_{0}^{t}\left\langle \frac{\partial}{\partial x}I_{\epsilon,1}(s,\cdot),I_{\epsilon,1}(s,\cdot)\right\rangle _{\tilde{L}_{\delta, w}^2}ds+\int_{0}^{t}\left\langle \frac{\partial}{\partial x}I_{\epsilon,h^{2}(z)}(s,\cdot),I_{\epsilon,1}(s,\cdot)\right\rangle _{\tilde{L}_{\delta, w}^2}ds \nonumber \\
& & \qquad -k_{1}\theta_{1}\int_{0}^{t}\left\langle \frac{\partial}{\partial y}I_{\epsilon, z^{-\frac{1}{2}}}(s,\cdot),I_{\epsilon,1}(s,\cdot)\right\rangle _{\tilde{L}_{\delta, w}^2}ds 
+k_{1}\int_{0}^{t}\left\langle \frac{\partial}{\partial y}I_{\epsilon,z^{\frac{1}{2}}}(s,\cdot),I_{\epsilon,1}(s,\cdot)\right\rangle _{\tilde{L}_{\delta, w}^2}ds \nonumber \\
& & \qquad +\int_{0}^{t}\left\langle \frac{\partial^{2}}{\partial x^{2}}I_{\epsilon,h^{2}(z)}(s,\cdot),\,I_{\epsilon,1}(s,\cdot)\right\rangle _{\tilde{L}_{\delta, w}^2}ds 
+ \frac{\xi_{1}^{2}}{4}\int_{0}^{t}\left\langle \frac{\partial^{2}}{\partial y^{2}}I_{\epsilon,1}(s,\cdot),
I_{\epsilon,1}(s,\cdot)\right\rangle _{\tilde{L}_{\delta, w}^2}ds \nonumber \\
& & \qquad +\rho\int_{0}^{t}\left\langle \frac{\partial^{2}}{\partial x \partial y}I_{\epsilon,h\left(z\right)}(s,\cdot),\,I_{\epsilon,1}(s,\cdot)\right\rangle _{\tilde{L}_{\delta, w}^2}ds \nonumber \\
& & \qquad +\xi_{1}\rho_{3}\rho_{1,1}\rho_{2,1}\int_{0}^{t}\left\langle \frac{\partial}{\partial x}I_{\epsilon,h\left(z\right)}(s,\cdot),\,\frac{\partial}{\partial y}I_{\epsilon,1}(s,\cdot)\right\rangle _{\tilde{L}_{\delta, w}^2}ds \nonumber \\
& & \qquad + \frac{\xi_{1}^{2}}{4}\int_{0}^{t}\left\langle \frac{\partial}{\partial y}I_{\epsilon, z^{-\frac{1}{2}}}(s,\cdot),I_{\epsilon,1}(s,\,\cdot)\right\rangle _{\tilde{L}_{\delta, w}^2}ds+\rho_{1,1}^{2}\int_{0}^{t}\left\Vert \frac{\partial}{\partial x}
I_{\epsilon,h(z)}(s,\cdot)\right\Vert _{\tilde{L}_{\delta, w}^2}^{2}ds \nonumber \\
& & \qquad + \frac{\xi_{1}^{2}}{4}\rho_{2,1}^{2}\int_{0}^{t}\left\Vert \frac{\partial}{\partial y}
I_{\epsilon,1}(s,\cdot)\right\Vert _{\tilde{L}_{\delta, w}^2}^{2}ds
-2\rho_{1,1}\int_{0}^{t}\left\langle \frac{\partial}{\partial x}I_{\epsilon,h(z)}(s,\cdot),I_{\epsilon,1}(s,\cdot)\right\rangle _{\tilde{L}_{\delta, w}^2}dW_{s}^{0} \nonumber \\
& &\qquad -\xi_{1}\rho_{2,1}\int_{0}^{t}\left\langle \frac{\partial}{\partial y}I_{\epsilon,1}(s,\cdot),I_{\epsilon,1}(s,\cdot)\right\rangle _{\tilde{L}_{\delta, w}^2}dB_{s}^{0}. \label{eq:5.11}
\end{eqnarray}
Observe now that by the definition of $u_{xx}$ in our SPDE, we have
\begin{eqnarray}
& & \int_{\mathbb{R}^+}\int_{\mathbb{R}}u_{xx}(s,\,x,\,z)\phi_{\epsilon}(z,\,y)w^2(x)f(x)dzdx \nonumber \\
& & \qquad\qquad = \int_{\mathbb{R}^+}\int_{\mathbb{R}}u(s,\,x,\,z)\phi_{\epsilon}(z,\,y)\left(w^2(x)f(x)\right)_{xx}dzdx \nonumber \\
& & \qquad\qquad =-\int_{\mathbb{R}^+}\int_{\mathbb{R}}u_{x}(s,\,x,\,z)\phi_{\epsilon}(z,\,y)\left(w^2(x)f(x)\right)_{x}dzdx \nonumber \\
& & \qquad\qquad =-\int_{\mathbb{R}^+}\int_{\mathbb{R}}w^2(x)u_{x}(s,\,x,\,z)\phi_{\epsilon}(z,\,y)f_{x}(x)dzdx \nonumber \\
& & \qquad\qquad \qquad -\int_{\left[0, \, 1\right]}\int_{\mathbb{R}}u_{x}(s,\,x,\,z)\phi_{\epsilon}(z,\,y)f(x)dzdx \label{eq:5.12}
\end{eqnarray}
which equals
\begin{eqnarray}
& & -\int_{\mathbb{R}^+}\int_{\mathbb{R}}w^2(x)u_{x}(s,\,x,\,z)\phi_{\epsilon}(z,\,y)f_{x}(x)dzdx \nonumber \\
& & \qquad +\int_{\left[0, \, 1\right]}\int_{\mathbb{R}}u(s,\,x,\,z)\phi_{\epsilon}(z,\,y)f_{x}(x)dzdx -\int_{\mathbb{R}}u(s,\,1,\,z)\phi_{\epsilon}(z,\,y)f(1)dz \nonumber
\end{eqnarray}
for any smooth function $f$ defined on $\left[0, \, +\infty\right)$. Since $u\in H_{\alpha}$ and since $f(1)$ can be controlled by the $H_{\alpha}$ norm of $f$ (by using Morrey's inequality near $1$), \eqref{eq:5.12} defines a linear functional on the space of smooth functions $f$ (defined on $\left[0, \, +\infty\right)$) which is bounded under the topology of $H_{\alpha}$. Then, since those functions form a dense subspace of $H_{\alpha}$, we have that \eqref{eq:5.12} holds also for any $f \in H_{\alpha}$. In particular, for $f=I_{\epsilon,1}(s,\,\cdot,\,y)$, multiplying \eqref{eq:5.12} by $y^{\delta}$ and then integrating in $\left(y, t\right)$ over $\mathbb{R}^{+} \times \mathbb{R}^{+}$, we obtain
\begin{eqnarray}
&&\int_{0}^{t}\left\langle \frac{\partial^{2}}{\partial x^{2}}I_{\epsilon,h^{2}(z)}(s,\,\cdot),\,I_{\epsilon,1}(s,\,\cdot)\right\rangle _{\tilde{L}_{\delta, w}^2}ds \nonumber \\
&& \qquad =-\int_{0}^{t}\left\langle \frac{\partial}{\partial x}I_{\epsilon,h^{2}(z)}(s,\,\cdot),\,\frac{\partial}{\partial x}I_{\epsilon,1}(s,\,\cdot)\right\rangle _{\tilde{L}_{\delta, w}^2}ds. \nonumber \\
&& \qquad \qquad -\int_{0}^{t}\left\langle \frac{\partial}{\partial x}I_{\epsilon,h^{2}(z)}(s,\,\cdot),\,\mathbb{I}_{\left[0, \, 1\right] \times \mathbb{R}}(\cdot)I_{\epsilon,1}(s,\,\cdot)\right\rangle _{\tilde{L}_{\delta}^2}ds. \label{eq:5.13}
\end{eqnarray}
Next, identities (5.10)-(5.12) from \cite{HK17} hold also when $\tilde{L}_{\delta'}^2$ is replaced by $\tilde{L}_{\delta', w}^2$ for $\delta' \in \{\delta, \delta - 1, \delta - 2\}$, and their justification is identical. Substituting these and \eqref{eq:5.13} in \eqref{eq:5.11} we obtain the desired result.
\end{proof}

\begin{rem}
It is equation~\eqref{eq:5.13} which adds two extra terms in the $\delta$-identity.
\end{rem}

{\flushleft{\textbf{Acknowledgement}}\\[.1in]
The second author's work was supported financially by the United Kingdom Engineering and Physical Sciences Research Council {[}EP/L015811/1{]}}, and by the Foundation for Education and European Culture in Greece (founded by Nicos \& Lydia Tricha).

\appendix

\section{APPENDIX: A clarification on the proof of Theorem 4.1}

The last computation in that proof assumes that $A$ is almost surely a continuity set of $X_{t}$. To see this, consider the process $Y$ satisfying the same SDE and initial condition as $X$ but without the stopping condition at $0$, and observe that it is a Gaussian process given the path $W_{\cdot}^{0}$ and given $\mathcal{G}$, which implies that 
\begin{equation}
\mathbb{P}\left(X_{t}\in V\,|\,W_{\cdot}^{0},\,\mathcal{G}\right) \leq \mathbb{P}\left(Y_{t}\in V\,|\,W_{\cdot}^{0},\,\mathcal{G}\right) = 0
\end{equation}
for any Borel set $V\subset\mathbb{R}^{+}$ of zero Lebesgue measure.

\end{document}